# Computational and analytical studies of a new nonlocal phase-field crystal model in two dimensions

Qiang Du

*Department of Applied Physics and Applied Mathematics, and Data Science Institute,*
*Columbia University, New York, NY 10027, USA*
*qd2125@columbia.edu*

Kai Wang

*School of Mathematics and Statistics,*
*Central South University, Changsha, Hunan, China*
*kai.wang.kw@csu.edu.cn*

Jiang Yang[*]

*Department of Mathematics, Shenzhen International Center for Mathematics & National*
*Center for Applied Mathematics Shenzhen (NCAMS),*
*Southern University of Science and Technology, Shenzhen, China*
*yangj7@sustech.edu.cn*

A nonlocal phase-field crystal (NPFC) model is presented as a nonlocal counterpart of the local phase-field crystal (LPFC) model and a special case of the structural PFC (XPFC) derived from classical field theory for crystal growth and phase transition. The NPFC incorporates a finite range of spatial nonlocal interactions that can account for both repulsive and attractive effects. The specific form is data-driven and determined by a fitting to the materials structure factor, which can be much more accurate than the LPFC and previously proposed fractional variant. In particular, it is able to match the experimental data of the structure factor up to the second peak, an achievement not possible with other PFC variants studied in the literature. Both LPFC and fractional PFC (FPFC) are also shown to be distinct scaling limits of the NPFC, which reflects the generality. The advantage of NPFC in retaining material properties suggests that it may be more suitable for characterizing liquid-solid transition systems. Moreover, we study numerical discretizations using Fourier spectral methods, which are shown to be convergent and asymptotically compatible, making them robust numerical discretizations across different parameter ranges. Numerical experiments are given in the two-dimensional case to demonstrate the effectiveness of the NPFC in simulating crystal structures and grain boundaries.

---

[*]Corresponding author.

## 1. Introduction

Over the years, various continuum models have been developed to describe different aspects of crystal growth and liquid-solid transitions. One such model that has generated considerable interest is the phase-field crystal (PFC) model proposed by Elder et al. [24,25]. The PFC model is a conserved form of the non-conserved Swift–Hohenberg equation [52]. It introduces a periodic order parameter $\phi$ to represent the local-time-averaged atomic density field, with the associated dimensionless free energy

$$\mathcal{E}(\phi) = \frac{1}{2}\|\Delta\phi + \phi\|_2^2 + \int_\Omega F(\phi), \tag{1.1}$$

where $\Omega$ is the unit period domain of $\phi$. Thus, the PFC model can be viewed as an $H^{-1}$-gradient flow of energy (1.1) with a constant mobility:

$$\frac{\partial \phi}{\partial t} = \Delta\mu, \quad \text{where} \quad \mu = (\Delta + I)^2 \phi + F'(\phi). \tag{1.2}$$

Although initially proposed to model elasticity in crystal growth, the PFC models have been applied to various other fields, such as thin film growth, dendrite formation, single dislocation, alloy solidification, and spontaneous elastic interaction. They have been successfully applied to study crystals with 2D triangular and 3D BCC symmetries. More studies and applications of PFC models can be found in recent reviews [27,39,51] and references provided therein. To allow for effective simulations of other common metallic crystal structures, various extensions have been studied. These extensions include the use of higher-order derivatives in free energy formulations [15, 29, 31, 37, 59] and structural PFC (XPFC) models that adopt a nonlocal integral form in the free energy formulation [30, 42, 43]. A special type of fractional nonlocal interactions has also been studied in [2, 3].

In this work, we focus on a nonlocal analog of the original phase-field crystal model (1.2). For easy reference, the latter model (1.2) is referred to as the local phase field crystal (LPFC) model and (1.1) as the local energy. We first define the nonlocal free energy

$$\mathcal{E}_\delta(\phi) = \frac{1}{2}\|\mathcal{L}_\delta\phi + \phi\|_2^2 + \int_\Omega F(\phi), \tag{1.3}$$

obtained by replacing the Laplace operator $\Delta$ in (1.1) with the nonlocal Laplacian, also named nonlocal diffusion (ND) operator, $\mathcal{L}_\delta$ parameterized by a constant $\delta$:

$$\mathcal{L}_\delta\phi(\mathbf{x}) = \int_{\mathcal{B}_\delta(\mathbf{x})} \rho_\delta(|\mathbf{x} - \mathbf{x}'|)\left(\phi(\mathbf{x}') - \phi(\mathbf{x})\right) d\mathbf{x}'. \tag{1.4}$$

Here, $\delta$ is the horizon parameter measuring the range of interactions, $\mathcal{B}_\delta(\mathbf{x})$ is a $\delta$-neighborhood of the center $\mathbf{x}$, and the kernel $\rho_\delta(s)$ is a nonnegative radial-type function with a compact support in $[0, \delta]$. Moreover, the kernel $\rho_\delta$ is taken such that we recover $\mathcal{L}_0 = \Delta$ as the limit of $\mathcal{L}_\delta$ as $\delta \to 0$.

Then, by introducing the nonlocal chemical potential using the variational derivative of the energy with respect to $\phi$,

$$\mu^\delta = (\mathcal{L}_\delta + I)^2 \phi + F'(\phi), \tag{1.5}$$

a nonlocal phase-field crystal (NPFC) model can be formulated as

$$\frac{\partial \phi}{\partial t} = \tilde{\mathcal{L}}_\delta \mu^\delta = \tilde{\mathcal{L}}_\delta((\mathcal{L}_\delta + I)^2 \phi + F'(\phi)), \tag{1.6}$$

where $\tilde{\mathcal{L}}_\delta$ is a given self-adjoint and positive-definite linear operator, NPFC can be mathematically interpreted as a gradient flow of the energy (1.3) associated with the inner product $(v, (-\tilde{\mathcal{L}}_\delta)^{-1} v)$, similar to the original PFC. While the operator $\tilde{\mathcal{L}}_\delta$ can take on a general form, it is assumed to satisfy, for any constant function $\phi = c$, $\tilde{\mathcal{L}}_\delta(\phi) = 0$, which in particular implies that $(\tilde{\mathcal{L}}_\delta f, 1) = (f, \tilde{\mathcal{L}}_\delta 1) = 0$ for any periodic function $f$, where $(\cdot, \cdot)$ is the standard $L^2$ inner product. This leads to the conservation of total mass,

$$\frac{d}{dt} \int_\Omega \phi(x, t) dx = 0, \quad \forall t \geq 0.$$

Thus, (1.6) represents conservative dynamics, just like its local counterpart (1.2). For computational convenience, we also assume that $\tilde{\mathcal{L}}_\delta$ is simultaneously diagonalizable with $\Delta$ and $\mathcal{L}_\delta$. Moreover, to have consistency with the local model, we assume further that $\mathcal{L}_0 = \Delta$ is also the limit of $\tilde{\mathcal{L}}_\delta$ as $\delta \to 0$. Some special cases of $\tilde{\mathcal{L}}_\delta$ include $\tilde{\mathcal{L}}_\delta = \Delta$ and $\tilde{\mathcal{L}}_\delta = \mathcal{L}_\delta$.

There are several motivations for the study of the NPFC model that involves the ND operator with a finite range of nonlocal interactions. Firstly, this choice is in line with the XPFC approach, which explores more general forms of free energy than the original LPFC models. Secondly, as the nonlocal interaction kernel becomes more localized, we can recover the local partial differential operators in the limit, which enables us to make connections between the NPFC model and the original PFC model. Thirdly, by choosing special fractional kernels and setting the interaction domain to be the entire space, we can also recover fractional PFC (FPFC) models as another limiting case. In addition, it is worth noting that nonlocal operators are integral operators that have been widely used to model problems in physics, chemistry, and materials science. Although the PFC model was originally derived from a phenomenological approach, it can be reinterpreted as a simplified and approximated version of the dynamic density functional theory (DDFT). In this sense, the NPFC model is in the same spirit while offering us more flexibility in choosing the interaction kernels. Moreover, the ND operators are of the convolution type and are diagonalizable in the Fourier space. Thus, fast algorithms proposed in [21] can be used to evaluate the Fourier symbols associated with various choices of kernels. This significantly reduces the computational complexities in the simulations of NPFC and makes the simulation costs comparable to that of simulating the LPFC models. Furthermore, one of the features of ND operators is that they avoid the explicit use of spatial derivatives, allowing us to handle more singular solutions. Unlike diffusive interfaces generated by LPFC models, the NPFC can capture sharp interfaces among bulk phases without causing further complications. In some cases, these sharp interfaces have important physical meanings.

The use of ND operators involving a finite range of interactions, as formulated by (1.4) has been widely studied in various fields [16]. These operators have been used to model nonlocal heat conduction [10], phase transitions [9,13,28], kinetic equations [34], nonlocal Dirichlet forms [7], and peridynamics [8,38,47,49,50], among others. Significant progress has been made in the rigorous mathematical analysis [4,5,12,17,18,26] and the algorithmic development of these operators [11,14,32,35,41,48,57,61]. In particular, asymptotically compatible (AC) schemes have been developed to preserve the limiting behavior of ND operators in discretizations [53,54]. [20] successfully extended this AC concept to semidiscrete Fourier spectral methods for solving nonlocal nonlinear Allen–Cahn equations. Moreover, [20] has derived a uniform and optimal error estimate of $O(\delta^2)$ for the convergence of numerical nonlocal solutions to the numerical solution of the local limit. However, the study of [20] was limited to one-dimensional problems.

In this paper, we carefully investigate NPFC models from both modeling and numerical simulation perspectives. We first employ the data-driven modeling approach to illustrate that a more accurate fitting of the structure factors of materials in the liquid state can be obtained by selecting special forms of the nonlocal interaction kernels used in the NPFC that involve both repulsive and attractive interactions. This desirable feature helps to demonstrate the advantages of NPFC the LPFC and fractional PFC models on the physical ground. We then develop suitable numerical algorithms for the efficient simulations of NPFC. For time discretization, we show that SAV schemes can be used to preserve the features and structures of the original continuum models. For spatial discretization, we utilize the Fourier spectral discretization and fast algorithms for the evaluation of nonlocal operators in the spectral space. Moreover, we prove that Fourier spectral methods are not only convergent but also give asymptotically compatible discretizations of the 2D NPFC models together with a uniform and optimal error estimate. Additionally, we present simulations with sharp interfaces by NPFC models with integrable kernels which allows for more singular solutions to be captured, demonstrating the capability of NPFC models in generating sharp interfaces instead of typical diffusive interfaces associated with LPFC models.

## 2. Physical motivation and model development

Let us first recall some basic features of the free energy functional and the material structure factor to motivate the development of NPFC, particularly with sign-changing kernels that account for both repulsive and attractive interactions.

### 2.1. *Free energy and structure factor*

During the transition from liquid to solid phases, there is a significant change in density. Specifically, the density is relatively uniform in the liquid phase but becomes spatially periodic in the solid phase. To model this behavior, a free energy functional given by Eq. (2.1) has been developed in [24, Sec.I(C)]:

$$\mathcal{E}(\phi) = \int_\Omega \left\{ \frac{\phi}{2} G(\Delta)\phi + F(\phi) \right\}. \tag{2.1}$$

The first term in the above energy functional is associated with the Laplacian operator $\Delta$. It induces spatial dependence and is responsible for the periodicity observed in the solid phase. The second term $\int_\Omega F(\phi)$ captures the thermodynamic properties of the system, such as the energy required for the formation of the solid-liquid interface and density changes. The model is also prescribed by specific forms of $G(\Delta)$, which can be determined by experimental data or theoretical considerations. Overall, the free energy functional in Eq. (2.1) yields a simple yet effective model for studying the liquid-solid transition in a system and has been widely used in the literature.

For elastic materials, the simplest possible forms of the free energy (2.1) capable of producing periodic structures have been constructed in [24, Sec.I(C)] via the following representation $G$:

$$G(\Delta) = \lambda(k_c^2 + \Delta)^2. \tag{2.2}$$

The parameters used in (2.2) are to fit the structure factor for $^{36}Ar$, which was determined by experiments at $85K$ as reported in [60]. That is, the fitting was done with respect to the wave number $k$ for the structure factor function

$$S(k) = \frac{1}{\hat{\omega}(k)}$$

where $\hat{\omega}(k) = G(-k^2) - \epsilon$ with

$$\epsilon := -F''(0)$$

being a constant with its physical significance proportional to the undercooling, i.e., $\epsilon \sim T_e - T$, where $T$ is the temperature of the system and $T_e$ is the equilibrium temperature at which the phase transition occurs. Here, the term $G(-k^2)$ that resulted from the Fourier transform of (2.2) is of the form

$$G(-k^2) = \lambda(k_c^2 - k^2)^2.$$

It is worth noting that the fitted curve was often able to match the experimental data up to the first peak, which has been the usual practice since small wave numbers carry more information on the material. For comparison, we briefly introduce the FPFC which is demonstrated in [3]. In the FPFC model, instead of Eq.(2.2), the following operator

$$G(\Delta, \gamma) = \lambda(k_c^2 + \Delta)^{2\gamma} \text{ or } G(\Delta, \gamma) = \lambda(k_c^{2\gamma} - (-\Delta)^\gamma)^2 \tag{2.3}$$

is considered with a parameter $\gamma$. Obviously, FPFC reduces to the classical PFC with $\gamma = 1$. Instead of using local operators in the above models, we can replace them with nonlocal operators defined in (1.4), that is,

$$G(\mathcal{L}_\delta) = \lambda(k_c^2 - \mathcal{L}_\delta)^2 \tag{2.4}$$

In particular, with kernels $\rho_\delta$ that are allowed to be sign-changing, we can incorporate both repulsive and attractive interactions and offer better approximations to the structure factor, thus making it possible to produce and match beyond a single peak. The latter has not been achieved in the literature before.

### 2.2. NPFC models with sign-changing kernels

To describe the 2D NPFC models with the relevant ND operator $\mathcal{L}_\delta$ defined on the domain $[0, L_x] \times [0, L_y]$, we first introduce a couple of truncated and normalized fractional kernels

$$\rho_{\alpha_j,\delta_j}(s) = \frac{2(4-\alpha_j)}{\pi \delta^{4-\alpha_j} |s|^{\alpha_j}} \chi_{[-\delta_j,\delta_j]}, \quad s \in (0,\delta_j),\ \delta_j \in (0,\infty),\ \alpha_j \in [0,4). \tag{2.5}$$

Now, we focus on the ND operator $\mathcal{L}_\delta$ involving a kernel represented by a linear combination of two fractional kernels with fractional powers $\alpha_1$ and $\alpha_2$ and horizon parameters $\delta_1$ and $\delta_2$, that is, the kernel $\rho_\delta$ in (1.4) is taken to be of the form

$$\begin{aligned}
\rho_{\alpha,\delta}(s) &= c_1 \rho_{\alpha_1,\delta_1}(s) - c_2 \rho_{\alpha_2,\delta_2}(s) \\
&= \frac{2(4-\alpha_1)c_1}{\pi} \frac{1}{\delta_1^{4-\alpha_1} s^{\alpha_1}} - \frac{2(4-\alpha_2)c_2}{\pi} \frac{1}{\delta_2^{4-\alpha_2} s^{\alpha_2}}, \quad s \in (0,\delta),
\end{aligned} \tag{2.6}$$

where $\alpha = \{\alpha_1, \alpha_2\}$ and $\delta = \min\{\delta_1, \delta_2\}$. To satisfy the normalized moment condition, we require that

$$c_1 - c_2 = 1.$$

It is trivial to check that for a pair of integers $(k, l)$, $e^{i\left(\frac{2\pi k}{L_x} x + \frac{2\pi l}{L_y} y\right)}$ is an eigenfunction of $\mathcal{L}_\delta$ with periodic boundary conditions. The corresponding eigenvalue is $\lambda_\delta(k, l) = c_1 \lambda_{\alpha_1,\delta_1}(k, l) - c_2 \lambda_{\alpha_2,\delta_2}(k, l)$ where

$$\lambda_{\alpha_j,\delta_j}(k,l) = \frac{2(4-\alpha_j)}{\pi \delta_j^2} \int_0^1 r^{1-\alpha_j} \int_0^{2\pi} \left( \cos\left(2\pi r \delta_j \sqrt{\frac{k^2}{L_x^2} + \frac{l^2}{L_y^2}} \cos\theta_j \right) - 1 \right) d\theta dr \tag{2.7}$$

for $j = 1, 2$. The evaluation of the above double integrals can be performed using the accurate and efficient hybrid algorithm proposed in [21]. For suitable parameter values, the kernel $\rho_{\alpha,\delta}$ could change sign, being positive around the origin but negative at $r = \delta$.

### 2.3. *Fitting the structure factor with NPFC having sign-changing kernels*

As an illustration of NPFC models, we investigate the fitting of the structure factor with sign-changing kernels. Comparisons with the fitting results associated with the local Laplace operator and the FPFC studied in [3] are also carried out. For the local operator $\mathcal{L}_0 = \Delta$, the parameters $\{\lambda, k_c\}$ are selected by fitting the functional form $G(\Delta)$ to the first-order peak in the experimental measurements [25] of the structure factor represented by

$$S(k) = \frac{1}{G(-k^2) - \epsilon}.$$

For the nonlocal operator $\mathcal{L}_\delta$, we are able to select parameters $\{\lambda, k_c, \alpha_1, \alpha_2, \delta_1, \delta_2, c_1, c_2\}$ by fitting the functional form to the first-order peak in the experimental measurements with sign-changing kernels, i.e., we fit

$$S(\lambda_{\alpha,\delta}) = \frac{1}{G(\lambda_{\alpha,\delta}) - \epsilon},$$

where parameters $k_c$ and $\epsilon$ are chosen in advance by referring to experimental data.

Similarly to earlier studies, we utilize experimental data on the liquid structure factor of $^{35}Ar$ at a temperature of 85K, which is near the melting point of argon, where the material remains in a liquid state. The fitted structure factor profiles are shown in Fig.1, where the experimental data is displayed with black dotted line. We observe that the fitted curves with a local operator (presented with a magenta dotted line) and FPFC (presented with blue solid-dotted line) match well with the experimental data up to the first peak. However, the fitted profile with sign-changing kernels (presented with a solid red line) agrees almost up to the first two peaks. In contrast, it should be noted that the fitted curves with a local operator and FPFC can only have one peak. Moreover, for small wave numbers $k$, the fitting by the nonlocal operator with sign-changing kernels, like that for the FPFC, is much closer to the experimental data than the case with a local operator. This demonstrates that the nonlocal operator with sign-changing kernels has an advantage in capturing the energy compared to the local operator. The improvement in fit quality is remarkable. As an illustration, we present here the results of the fitting and the values of parameters used when $k_c = 1.997, \epsilon = -0.3725$: for the local case we have $\lambda = 0.8115$ as the only fitting parameter; for the FPFC case, the two parameters are $\lambda = 0.8666$ and $\gamma = 1.237$; meanwhile, for the nonlocal case, the parameters used are, respectively, $\lambda = 1.8315, \delta_1 = 2.21, \alpha_1 = 1, \delta_2 = 3.01, \alpha_2 = 0, c_1 = 3$, and $c_2 = 2$. The combined kernel with parameters resulted from fitting is checked to be positive, and the resulted operator is still negative definite.

The corresponding structure factors fitting of some other materials reported in [23, 58] are presented in Fig.2. In each case, the improvement made by the NPFC model with sign-changing kernels is striking compared with the LPFC and FPFC models. This shows that the resulting NPFC model will have better performance in capturing phases in modeling liquid-solid transition systems.

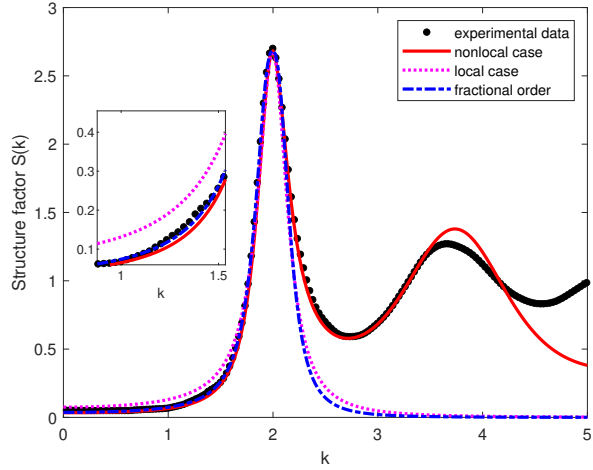

Fig. 1. Fit of structure factor of $^{36}Ar$.

Next, we make some comparisons with XPFC models with multi-peak Gaussians. For XPFC models [6], the free energy is expressed as

$$\mathcal{E}(\phi) = \int_\Omega \left[ \frac{1-\varepsilon}{2}\phi^2 - \frac{a}{3}\phi^3 + \frac{v}{4}\phi^4 \right] - \frac{1}{2}\iint_\Omega [\phi(x)C_2(\|x-x'\|)\phi(x')]\,dx'dx,$$

where the pair correlation function is approximated by a combination of modulated Gaussian functions in Fourier space via

$$\hat{C}_2(k) = \max(G^i(k), G^{i+1}(k), \ldots, G^N(k)),$$

with $N$ being the total number of Gaussian functions used in the approximation of the direct correlation function, and

$$G^i(k) = \exp\left(-\frac{\sigma^2 k_i^2}{2\lambda_i \beta_i}\right) \exp\left(-\frac{(k-k_i)^2}{2\alpha_i^2}\right)$$

the modulated Gaussian function (i.e., a Gaussian function with its height modified by an exponential function). The parameter $k_i$ specifies the position of the $i^{\text{th}}$ Gaussian peak, $\alpha_i$ corresponds to the root-mean-square width of the $i^{\text{th}}$ Gaussian peak and controls the excess energy associated with defects, interfaces, and strain, $\sigma$ controls the heights of the Gaussian peaks and is related to temperature, $\lambda_i$ and $\beta_i$ are the planar atomic density and the number of planar symmetries of the $i^{\text{th}}$ family of crystallographic planes. In this XPFC model, the structure factor $S(k)$ can be obtained from the peak position via

$$S(k) = \frac{1}{N}\sum_{i=1}^N \sigma_i \exp\left(-\frac{(k-k_i)^2}{2\alpha_i^2}\right).$$

Note that here we simplify the parameters to $\sigma_i$ associated with the height of the peak and $\alpha_i$ associated with the width of the peak. The results are shown in Fig.3.

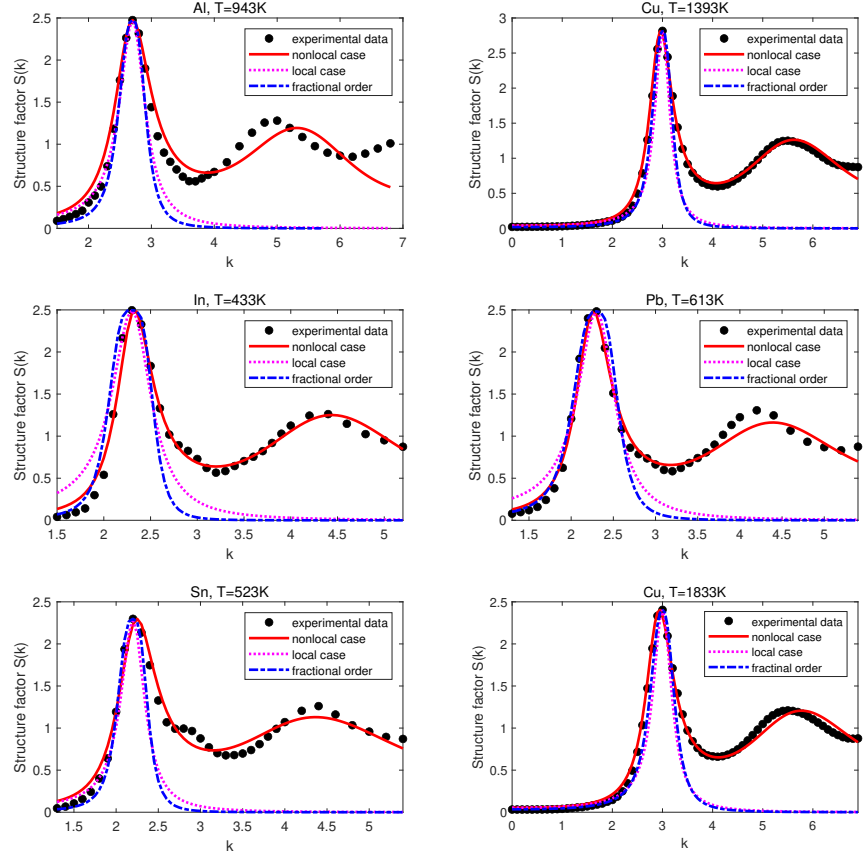

Fig. 2. Structure fitting for different materials at various temperatures.

Here we set $N = 2$ to fit the first two peaks with $k_1 = 1.9971$ and $k_2 = 3.6417$. While XPFC models using multi-peak Gaussians can fit the structural factor with more peaks, NPFC models not only offer a much better fitting than the XPFC around the second peak but also lead to slightly better fitting even at the first peak. This provides further support for the use of NPFC over XPFC variants of the PFC.

Note that linear combinations of truncated power-like kernels are adopted for nonlocal operators. They are specifically chosen to achieve good fittings of structure factors. Despite the sign changes of the kernels, the resulting nonlocal operators can remain positive definite; see [36] for related discussions. One may verify this easily for the specific choices used here by evaluating their eigenvalues that can be efficiently obtained, see later discussions. In the rest of this paper, without ambiguity, we always simplify $\rho_{\alpha_j,\delta_j}$ to $\rho_\delta$.

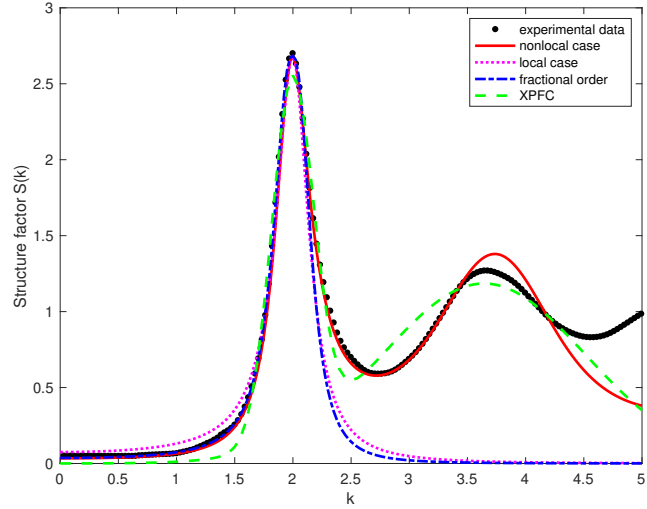

Fig. 3. Fit of structure factor of $^{36}Ar$ with XPFC.

## 3. Properties of the NPFC and the discretization schemes

To present some theoretical analysis, we assume that the operators $-\mathcal{L}_\delta$ and $-\tilde{\mathcal{L}}_\delta$ are all positive definite, which are reasonable assumptions that can be verified for the fitting parameters used.

The NPFC equation given in (1.6) is a gradient flow of nonlocal free energy (1.3) associated with the dual space norm of the energy norm associated with the operator $-\tilde{\mathcal{L}}_\delta$. By the assumption on its positive-definiteness, we have the following energy identity

$$\frac{d\mathcal{E}_\delta(\phi^\delta)}{dt} = (\tilde{\mathcal{L}}_\delta \mu^\delta, \mu^\delta) \leq 0, \tag{3.1}$$

which is analogous to the LPFC equation, an $H^{-1}$ gradient flow satisfying

$$\frac{d\mathcal{E}(\phi^0)}{dt} = (\Delta \mu^0, \mu^0) \leq 0. \tag{3.2}$$

We note that the nonlocal chemical potential $\mu^\delta$ is given by (1.5) and $\mu^0$ is the local chemical potential defined analogously with $\Delta$ in place of $\mathcal{L}_\delta$. They are also the first variations of the energy functionals (1.3) and (1.1) respectively. If we define

$$\mathcal{G}_0 \phi = (\Delta + I)^2 \phi, \quad \text{and} \quad \mathcal{G}_\delta \phi^\delta = (\mathcal{L}_\delta + I)^2 \phi^\delta,$$

the NPFC (1.6) and the associated LPFC model can be respectively written as

$$\frac{\partial \phi^\delta}{\partial t} = \tilde{\mathcal{L}}_\delta(\mathcal{G}_\delta \phi^\delta + F'(\phi^\delta)), \quad \text{and} \quad \frac{\partial \phi^0}{\partial t} = \Delta(\mathcal{G}_0 \phi^0 + F'(\phi^0)). \tag{3.3}$$

Here, we discuss the nonlinear term $F(\phi)$ in more detail. The energy dissipation law of local models has sufficient regularities to guarantee the $L^\infty$-boundedness of both the exact solution and numerical solutions due to the Sobolev embedding $H^2 \hookrightarrow L^\infty$, which plays a crucial role in proving the error bounds. However, for nonlocal models, the energy dissipation does not have sufficient regularities to derive such boundedness. In order to derive the $L^\infty$-boundedness for nonlocal models following the Sobolev embedding $H^2 \hookrightarrow L^\infty$, we need to handle the term $\Delta F'(\phi)$. Therefore, we make the following assumption on the nonlinear term:

**Assumption** $$F \in C^2(\Omega) \quad \text{and} \quad \|F'''\|_\infty < C. \tag{3.4}$$

However, for the quartic polynomials, such as the typical choice $F(\phi) = \frac{1}{4}(\phi^2 - \epsilon)^2$, the assumption is not satisfied. To address this, we perform a Lipschitz truncation modification. Choosing a sufficiently large number $D$, we modify $F(\phi)$ with $\tilde{F}(\phi)$ as follows:

$$\tilde{F}(\phi) = \begin{cases} \frac{1}{4}(\phi^2 - \epsilon)^2, & |\phi| \leq D, \\ a\phi^2 + b \log|\phi| + ce^{-\phi^2} & |\phi| > D, \end{cases} \tag{3.5}$$

where coefficients $a, b$ and $c$ are determined such that $\tilde{F}(\phi)$, $\tilde{F}'(\phi)$ and $\tilde{F}''(\phi)$ are continuous at $|\phi| = D$. It is easy to check that with such a modification $\tilde{F}(\phi)$ (still denoted as $F(\phi)$ subsequently) satisfies the assumption (3.4).

We remark that in our numerical tests we set $D = 1$ and in fact, the maximum value of the numerical solution is always smaller than $D$. Hence, this modification does not actually affect the numerical simulation.

### 3.1. *Fourier spectral approximation in space*

We solve both LPFC and NPFC models numerically using Fourier spectral methods. The numerical solutions are of the following form

$$\phi_{MN}^\delta(x,y) = \sum_{|k| \leq M} \sum_{|\ell| \leq N} a_{k\ell}^\delta(t) e^{i\left(\frac{2\pi k}{L_x} x + \frac{2\pi \ell}{L_y} y\right)},$$

$$\phi_{MN}^0(x,y) = \sum_{|k| \leq M} \sum_{|\ell| \leq N} a_{k\ell}^0(t) e^{i\left(\frac{2\pi k}{L_x} x + \frac{2\pi \ell}{L_y} y\right)}$$

and satisfy

$$\begin{aligned} \frac{\partial \phi_{MN}^\delta}{\partial t} &= \tilde{\mathcal{L}}_\delta \mathcal{G}_\delta \phi_{MN}^\delta + \Pi_h[\tilde{\mathcal{L}}_\delta F'(\phi_{MN}^\delta)], \\ \frac{\partial \phi_{MN}^0}{\partial t} &= \mathcal{L}_0 \mathcal{G}_0 \phi_{MN}^0 + \Pi_h[\mathcal{L}_0 F'(\phi_{MN}^0)] \end{aligned} \tag{3.6}$$

respectively, with same initial value $\phi_0$, where $\Pi_h$ is the spectral orthogonal projection onto

$$S_{MN} = \text{span}\left\{ e^{i\left(\frac{2\pi k}{L_x} x + \frac{2\pi \ell}{L_y} y\right)} \,\bigg|\, -M \leq k \leq M, -N \leq \ell \leq N \right\},$$

defined for any $\phi$ by $\Pi_h \phi \in S_{MN}$ that satisfies

$$(\Pi_h \phi, \varphi) = (\phi, \varphi), \quad \forall \varphi \in S_{MN}. \tag{3.7}$$

Note that $\Pi_h$ commutes with $\tilde{\mathcal{L}}_\delta$ and $\mathcal{L}_0 = \Delta$.

### 3.2. *Fully discrete SAV scheme for time discretization*

We now discuss the time discretization of the time-dependent NPFC model. First, the NPFC models are reformulated with some scalar auxiliary variables (SAV), which follow steps proposed in [45, 46]:

$$\begin{cases} \dfrac{\partial \phi}{\partial t} = \tilde{\mathcal{L}}_\delta \mu, \\ \quad \mu = \mathcal{G}_\delta \phi + \beta \phi + r H(\phi), \\ \dfrac{dr}{dt} = \dfrac{1}{2} \int_\Omega H(\phi) \dfrac{\partial \phi}{\partial t} dx, \end{cases} \quad (3.8)$$

with

$$H(\phi) = \dfrac{U(\phi)}{\sqrt{\int_\Omega F(\phi)\mathrm{d}x - \frac{\beta}{2}\|\phi\|_2^2 + C_H}} \text{ and } U(\phi) = F'(\phi) - \beta\phi, \quad (3.9)$$

where $C_H$ and $\beta$ are any chosen positive constants.

The energy dissipation law follows by noting that

$$\dfrac{d\mathcal{E}(\phi(t))}{dt} = (\mu, \tilde{\mathcal{L}}_\delta \mu) \leq 0, \text{ with } \mathcal{E}(\phi) = \dfrac{1}{2}(\phi, \mathcal{G}_\delta \phi) + \dfrac{\beta}{2}(\phi, \phi) + r^2,$$

which is derived by taking inner products of the first two equations in scheme (3.8) against $\mu, \frac{\partial \phi}{\partial t}$, respectively, multiplying the third equation with $2r$ and summing up the resulted equations.

By making discretization with second order backward difference formula (BDF2) in time and Fourier spectral method in space, the fully discrete SAV/BDF2 scheme is given in (3.10) below:

$$\begin{cases} (3\phi_{MN}^{n+1} - 4\phi_{MN}^n + \phi_{MN}^{n-1}, q) = 2\tau(\tilde{\mathcal{L}}_\delta \mu_{MN}^{n+1}, q), \forall q \in S_{MN}, \\ (\mu_{MN}^{n+1}, w) = (\mathcal{G}_\delta \phi_{MN}^{n+1}, w) + \beta(\phi_{MN}^{n+1}, w) + r^{n+1}\left(H(\bar{\phi}_{MN}^{n+1}), w\right), \forall w \in S_{MN} \\ 3r^{n+1} - 4r^n + r^{n-1} = \dfrac{1}{2}\left(H(\bar{\phi}_{MN}^{n+1}), 3\phi_{MN}^{n+1} - 4\phi_{MN}^n + \phi_{MN}^{n-1}\right), \end{cases} \quad (3.10)$$

where

$$\bar{\phi}_{MN}^{n+1} = 2\phi_{MN}^n - \phi_{MN}^{n-1}. \quad (3.11)$$

### 3.3. *Energy stability and error estimates for the SAV/BDF2 scheme*

Multiplying the above three equations with $\mu_{MN}^{n+1}, (3\phi_{MN}^{n+1} - 4\phi_{MN}^n + \phi_{MN}^{n-1})/(2\tau)$ and $r^{n+1}/\tau$, respectively, integrating the first two equations and using the following elementary identity

$$\begin{aligned} 2(a^{k+1}, 3a^{k+1} - 4a^k + a^{k-1}) & \quad (3.12) \\ = |a^{k+1}|^2 + |2a^{k+1} - a^k|^2 + |a^{k+1} - 2a^k + a^{k-1}|^2 - |a^k|^2 - |2a^k - a^{k-1}|^2, \end{aligned}$$

we obtain

$$\frac{1}{\tau}\left\{\tilde{\mathcal{E}}[(\phi_{MN}^{n+1}, r^{n+1}), (\phi_{MN}^n, r^n)] - \tilde{\mathcal{E}}[(\phi_{MN}^n, r^n), (\phi_{MN}^{n-1}, r^{n-1})]\right\}$$
$$+ \frac{1}{\tau}\left\{\frac{1}{4}\left(\phi_{MN}^{n+1} - 2\phi_{MN}^n + \phi_{MN}^{n-1}, \mathcal{G}_\delta(\phi_{MN}^{n+1} - 2\phi_{MN}^n + \phi_{MN}^{n-1})\right)\right\}$$
$$+ \frac{1}{\tau}\left\{\frac{1}{4}\left(\phi_{MN}^{n+1} - 2\phi_{MN}^n + \phi_{MN}^{n-1}, \beta(\phi_{MN}^{n+1} - 2\phi_{MN}^n + \phi_{MN}^{n-1})\right)\right\}$$
$$+ \frac{1}{\tau}\left\{\frac{1}{2}(r^{n+1} - 2r^n + r^{n-1})^2\right\} = (\mu_{MN}^{n+1}, \tilde{\mathcal{L}}_\delta \mu_{MN}^{n+1}),$$

where for the positive semidefinite operator $\mathcal{G}_\delta$, the modified discrete energy is defined as

$$\tilde{\mathcal{E}}[(\phi_{MN}^{n+1}, r^{n+1}), (\phi_{MN}^n, r^n)] = \frac{1}{2}\left((r^{n+1})^2 + (2r^{n+1} - r^n)^2\right)$$
$$+ \frac{1}{4}\left((\phi_{MN}^{n+1}, \mathcal{G}_\delta \phi_{MN}^{n+1}) + (2\phi_{MN}^{n+1} - \phi_{MN}^n, \mathcal{G}_\delta(2\phi_{MN}^{n+1} - \phi_{MN}^n))\right)$$
$$+ \frac{1}{4}\left(\beta(\phi_{MN}^{n+1}, \phi_{MN}^{n+1}) + \beta(2\phi_{MN}^{n+1} - \phi_{MN}^n, 2\phi_{MN}^{n+1} - \phi_{MN}^n)\right),$$

which is always non-negative. Thus, with $-\tilde{\mathcal{L}}_\delta$ being positive definite, we get the following theorem.

**Theorem 3.1.** *The SAV/BDF2 scheme* (3.10) *is unconditionally energy stable in the sense that the modified energy* $\tilde{\mathcal{E}}[(\phi_{MN}^{n+1}, r^{n+1}), (\phi_{MN}^n, r^n)]$ *is non-increasing in* $n$.

In what follows we always assume $M = N$ for simplicity and use notation $u_N$ instead of $u_{MN}$ to represent the spatial discrete solution for any variable $u$. Without specific statement $\|\cdot\|$ stands for $L^2$ norm hereafter. We use $C$ to denote a generic positive constant that may take different values in its each occurrence, but is always independent of the temporal step size and spatial step size.

**Lemma 3.1.** *Let* $\{\phi_N^n, r^n\}$ *be solutions to scheme* (3.10) *with* $\tilde{\mathcal{L}}_\delta = \mathcal{L}_\delta$ *and let the assumption* (3.4) *hold. Assuming that* $\phi_0 \in H^2(\Omega)$, *then numerical solutions* $\{\phi_N^n, r^n\}$ *are bounded, that is,*

$$\|\phi_N^n\|_{L^\infty} + |r^n| \leq C. \tag{3.13}$$

The proof of the lemma is left to the Appendix 5.1. We now present a technical lemma to handle the operation of the nonlocal operators on nonlinear terms, which will be used to prove Theorem 3.2.

**Lemma 3.2.** *Assume that the functions $u$ and $v$ satisfies $\|u\|_\infty \leq C$ and $\|v\|_\infty \leq C$ and either $u$ or $v$ is globally Lipschitz continuous. Let $g = g(x)$ be any Lipschitz continuous function in $x \in [-C, C]$ and its derivative $g'$ is also Lipschitz continuous.*

*Define $e = u - v$, and then the following inequality holds*

$$(-\mathcal{L}_\delta(g(u) - g(v)), g(u) - g(v)) \leq C((-\mathcal{L}_\delta e, e) + (e, e)).$$

**Remark 3.1.** In [44], the following inequality, analogous to that in the Lemma 3.2,

$$\|\nabla g(u) - \nabla g(v)\| \leq C(\|e\| + \|\nabla e\|)$$

is proved with the assumptions $|g'(x)-g'(y)| \leq L|x-y|$ and $\|g\|+\|v\| \leq C$. However, the chain rule used in [44] is not applicable to the nonlocal operator considered here. Thus, we establish the conclusion in the above lemma subject to the assumption on the Lipschitz continuity of either $u$ or $v$, or its convex linear combination, see the proof in the appendix. In the application of the lemma, we take one of the $u$ and $v$ as the true solution for which the Lipschitz continuity can be assumed or rigorously established.

Next, we present an error estimate of the fully discrete SAV/BDF2 scheme for the NPFC model, under sufficient regularity assumptions on the true solution $\phi$.

**Theorem 3.2.** Let $\{\phi(t,x), r(t)\}$ and $\{\phi_N^n, r^n\}$ be solutions to Eqs.(3.8) and (3.10), respectively, with $\tilde{\mathcal{L}}_\delta = \mathcal{L}_\delta$. Suppose that $r \in H^3(0,T)$, $\phi_0 \in H^2(\Omega)$, $\partial_t^2 \phi \in L^2(0,T; H^1_{per}(\Omega))$, $\partial_t^3 \phi \in L^2(0,T; H^{-1}(\Omega))$, and the assumption (3.4) hold, then the following error estimation holds

$$\|\phi_N^n - \phi(t_n)\| + |r^n - r(t_n)| \leq C(\|\Pi_h \phi(t_n) - \phi(t_n)\| + \tau^2),$$

and in particular,

$$\|\phi_N^n - \phi(t_n)\| + |r^n - r(t_n)| \leq C(h^m + \tau^2),$$

for $\phi \in L^\infty(0,T; H^m_{per}(\Omega))$ and

$$\|\phi_N^n - \phi(t_n)\| + |r^n - r(t_n)| \leq C(e^{-c/h} + \tau^2),$$

for analytical solution $\phi$, where constants $C$ and $c$ are independent of the temporal step size $\tau$ and spatial step size $h$.

**Remark 3.2.** The proof of the above theorem is given in the Appendix 5.3. We note that the conclusion can be directly generalized to the case with any positive definite operator $-\tilde{\mathcal{L}}_\delta$ and non-negative definite operator $\mathcal{G}_\delta$. Moreover, an error estimate of $\mu$, omitted in the statement of the theorem, can also be obtained from the proof.

## 4. Asymptotic compatibility

It is known that the horizon parameter $\delta$ in the definition of nonlocal operator measures the range of nonlocal interactions. With $\delta \to 0$, only local interactions take effect, that is, the limit of nonlocal operator as $\delta$ goes to zero is the local differential operator. Numerical schemes that preserve this limit behavior are called asymptotically compatible schemes [53, 54]. Mathematically, let $u^0$ and $u_N^\delta$ be solutions to a continuous local model and its discrete nonlocal counterpart, respectively, then the asymptotic compatibility can be formulated as

$$\|u_N^\delta - u^0\| \to 0 \quad \text{as} \quad \delta \to 0, \, N \to \infty. \tag{4.1}$$

Further the estimate can be derived by following the triangle inequality below

$$\|u_N^\delta - u^0\| \leq \|u_N^\delta - u_N^0\| + \|u_N^0 - u^0\|. \tag{4.2}$$

The reason why we adopt such a triangle inequality is explained in [20].

In this section, we aim to prove the asymptotic compatibility of the discretization of the 2D NPFC models. The analysis is similar to that studied for nonlocal Allen-Cahn equations [20], so some details are skipped to keep the presentation short. Note that similar techniques can be used to establish the well-posedness of the NPFC model as well as its local limit as $\delta \to 0$.

### 4.1. *Analysis of asymptotic compatibility*

The asymptotic compatibility for the fully discrete scheme is presented in Theorem 4.2. Before presenting the theorem, we introduce Lemma 4.1 and Theorem 4.1 which will be used in the proof of Theorem 4.2.

**Lemma 4.1.** *Let $\phi_N^\delta(x,y)$ and $\phi_N^0(x,y)$ be numerical solutions to the steady linear nonlocal problem*

$$-\mathcal{L}_\delta \phi_\delta(x,y) = f(x,y) \quad \forall (x,y) \in (-\pi,\pi) \times (-\pi,\pi),$$

*and local problem*

$$-\mathcal{L}_0 \phi_0(x,y) = f(x,y) \quad \forall (x,y) \in (-\pi,\pi) \times (-\pi,\pi),$$

*respectively, subject to periodic conditions. Thus, it holds that*

$$\|\phi_N^\delta - \phi_N^0\| \leq C\delta^2 \|f\|,$$

*where the generic constant $C$ is independent of spatial step size.*

The result can be seen as an extension of similar conclusions presented in the one-dimensional case in [20]. The proof can be found in the Appendix 5.4, together with the proof of the following theorem.

**Theorem 4.1.** *Assume that $\phi_N^\delta$ and $\phi_N^0$ are solutions to (3.6), respectively, and and the assumption (3.4) hold. Then for any finite $t \leq T$, it holds that*

$$\|\phi_N^\delta(t,\cdot) - \phi_N^0(t,\cdot)\|_2 \leq C(T,\phi_0)\delta^2, \tag{4.3}$$

*where $C$ is independent of spatial step size.*

**Theorem 4.2 (Asymptotic compatibility).** *Let $\{\phi_N^{\delta,n+1}, r^{\delta,n+1}\}$ be numerical solutions to the above fully discrete nonlocal SAV scheme (3.10) and $\{\phi_N^{0,n+1}, r^{0,n+1}\}$ to corresponding local scheme, respectively. Suppose $\phi_0 \in H^2(\Omega)$. Then the following asymptotic convergence holds*

$$\|\phi_N^{\delta,n+1} - \phi_N^{0,n+1}\| + |r^{\delta,n+1} - r^{0,n+1}| \leq C(T,\phi_0)\delta^2$$

*where the constant $C$ is independent of spatial step size and temporal step size.*

**Remark 4.1.** Theorems proved above can be generalized to any negative definite operator $\tilde{\mathcal{L}}_\delta$ and non-negative definite operator $\mathcal{G}_\delta$.

**Remark 4.2.** In proving asymptotic compatibility, the requirement of global Lipschitz continuity assumption on nonlinear term $F$ can be relaxed to local Lipschitz. In this setting, the locally Lipschitz continuous nonlinear term $F(\phi)$ can be modified with the cut-off function Eq.(3.5), and the resulted nonlinear term $\tilde{F}(\phi) = F(\phi)$ for considered domain $\phi \in D$. Thus all remaining details hold and so is the conclusion.

### 4.2. *Numerical implementation*

In the following numerical examples, we will use the free energy (1.3) and the residual of the equation to decide whether the evolution is close enough to the steady state or not. In all numerical examples presented in this paper, we take quartic polynomial form of the nonlinear term $F(\phi) = \frac{1}{4}(\phi^2 - \epsilon)^2$. We approximate the continuum energy (1.3) as

$$E_{MN} = \frac{L_x L_y}{2} \sum_{|m| \leq M} \sum_{|n| \leq N} |1 + \lambda_\delta(m,n)|^2 |a_{mn}|^2 + \frac{1}{4(2M+1)(2N+1)} \sum_{j=0}^{2M} \sum_{k=0}^{2N} |\phi_{jk}^2 - \epsilon|^2,$$

where

$$\phi_{jk} = \sum_{|m| \leq M} \sum_{|n| \leq N} a_{mn} e^{i\left(\frac{2\pi m}{L_x} x_j + \frac{2\pi n}{L_y} y_k\right)}, \quad x_j = \frac{jL_x}{2M+1}, \quad y_k = \frac{kL_y}{2N+1}.$$

The residual is calculated as

$$\text{Res} = \left\| \frac{3\Phi^{n+1} - 4\Phi^n + \Phi^{n-1}}{2\tau} \right\|_F,$$

where $\|\cdot\|_F$ denotes Frobenius norm and $\Phi = \{\phi_{jk}\}$ is a $(2M+1) \times (2N+1)$ matrix.

### 4.3. *Numerical tests on the asymptotic compatibility for NPFC models*

**Example 4.1.** We firstly verify the asymptotic compatibility for a 2D NPFC model in $[0, 10\pi] \times [0, 10\pi]$ with periodic boundary conditions beiing considered in this example and initial value is given by $\phi_0(x) = 0.05 \sin(x) \cos(x) + 0.07$. We fix $\epsilon = 0.025$ and take the kernel defined in (2.5) with $\alpha = 1$ and $\alpha = 3$ respectively.

In both cases of $\alpha = 1$ and $\alpha = 3$, we take a small time step with $\Delta t = 0.001$ and $M = N = 64$. The numerical differences between the solutions of the NPFC and LPFC models are shown in Tab.1. The convergence speed is seen to be about $O(\delta^2)$, which is in good agreement with the theoretical results. Thus, we achieve the optimal convergence rate for asymptotic compatibility.

Table 1. (Example 4.1) Errors between numerical solutions of NPFC and LPFC models.

|  | $\alpha = 1$ | | $\alpha = 3$ | |
| --- | --- | --- | --- | --- |
| $\delta = 0.4$ | $\|\phi_N^\delta - \phi_N^0\|_2$ | Rate | $\|\phi_N^\delta - \phi_N^0\|_2$ | Rate |
| $\delta$ | 3.49e-04 | - | 1.71e-04 | - |
| $\delta/2$ | 7.08e-05 | 2.30 | 3.81e-05 | 2.16 |
| $\delta/4$ | 1.68e-05 | 2.07 | 9.27e-06 | 2.04 |
| $\delta/8$ | 4.15e-06 | 2.01 | 2.30e-06 | 2.01 |

### 4.4. *Numerical experiments*

Next, we present examples that offer comparisons between lattices simulated with LPFC and NPFC models respectively.

**Example 4.2.** Consider both 2D NPFC and LPFC models in $[0, 50] \times [0, 50]$ with periodic boundary conditions and the random initial value (RIV) with average mean value $\bar{\phi} = 0.07$. In both nonlocal and local models we fix $\epsilon = 0.025$. The fractional power kernel defined in (2.5) is taken as $\alpha = 1$.

We set $\Delta t = 1$, $M = N = 64$. In the nonlocal cases we choose small $\delta = 0.2$ and large $\delta = 2$. The numerical solutions at $T = 250, 500, 1000, 2000$ are presented in Fig.4. It is observed that the dynamics and steady state of NPFC models with small $\delta$ look nearly identical to the LPFC model. However, with a large $\delta$, there are visible differences between the dynamics and steady state of the NPFC models and the local ones. However, without considering the slight rotation, the steady states all have the same periodic patterns.

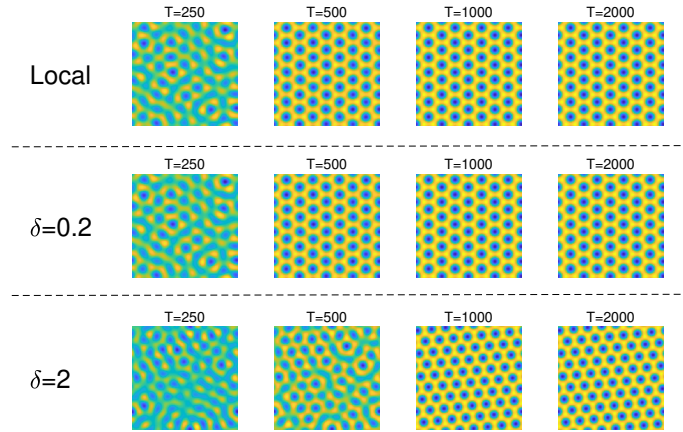

Fig. 4. (Example 4.2) Numerical evolutions of NPFC and LPFC models.

### 4.5. *Grain boundary simulation*

**Example 4.3.** Simulations of the grain boundaries with NPFC model are presented in this example. On a period grid of $[0, L_x] \times [0, L_y]$, we specify a one-mode approximation of the stationary hexagonal state $\phi_h(x, y)$ in one orientation between $0 < y < L_y/2$, while another hexagonal state of a different orientation is specified between $L_y/2 < y < L_y$. The hexagonal lattice is given by

$$\phi(x, y) = \bar{\phi} + A_\phi \left[ \cos(qx) \cos\left(\frac{qy}{\sqrt{3}}\right) - \frac{1}{2} \cos\left(\frac{2qy}{\sqrt{3}}\right) \right]$$

with $q = \frac{\sqrt{3}}{2}$, $\bar{\phi}$ being the mean value of $\phi$ and

$$A_\phi = \frac{4}{5} \left( \bar{\phi} + \frac{1}{3} \sqrt{15\epsilon - 36\bar{\phi}^2} \right).$$

In the middle of two hexagonal states with various orientations, random values with mean value $\bar{\phi}_l = \bar{\phi}$ is specified so as not to influence the nature of the grain boundary that emerged. The parameters for the simulation are taken as: $\epsilon = 4/15, \bar{\phi} = 0.2$, the time step size is 0.5 and the space step size is around $\pi/4$.

The results of numerical simulations corresponding to different values of grain mismatch angle $\theta = 11.6°, 26.3°, 39.4°$ at time $T = 10000$ are shown in Fig.5. It is observed that around interfaces 5|7 dislocation dipoles are formed instead of the original hexagonal lattice. We mark 5|7 dislocation dipoles with green pentagons and red heptagons in the corresponding figures for the ease of reading. These phenomena are similar to the results presented in [2].

### 4.6. *NPFC models with integrable kernels*

In [20], it was proved that NAC equations with integrable kernels admit discontinuous steady states under reasonable assumptions. Here we extend such a result to NPFC models.

First, the steady state of the conserved NPFC model with periodic boundary conditions satisfies

$$\mathcal{L}_\delta^2 \phi + 2\mathcal{L}_\delta \phi + (1 - \epsilon)\phi + \phi^3 = c_0, \tag{4.4}$$

where $c_0$ is a constant. On the other hand, for integrable kernels, we split the ND operator as

$$\mathcal{L}_\delta \phi = \rho_\delta \star \phi - c_\delta \phi, \tag{4.5}$$

where $c_\delta = \int_{\mathcal{B}_\delta(\mathbf{0})} \rho_\delta(|\boldsymbol{x}|) d\boldsymbol{x} > 0$, and $\star$ represents the convolution. Then, the steady state equation becomes

$$\phi^3 + \left((c_\delta - 1)^2 - \epsilon\right) \phi = -\rho_\delta \star (\rho_\delta \star \phi) + 2c_\delta (\rho_\delta \star \phi) - 2\rho_\delta \star \phi + c_0. \tag{4.6}$$

**Example 4.4.** In this example, we fix $\alpha = 0.5$ and $\delta = 3$ in the fractional kernel (2.5). The domain is taken as $[0, 12\pi] \times [0, 12\pi]$, and the discretization parameters are taken as $M = N = 128$.

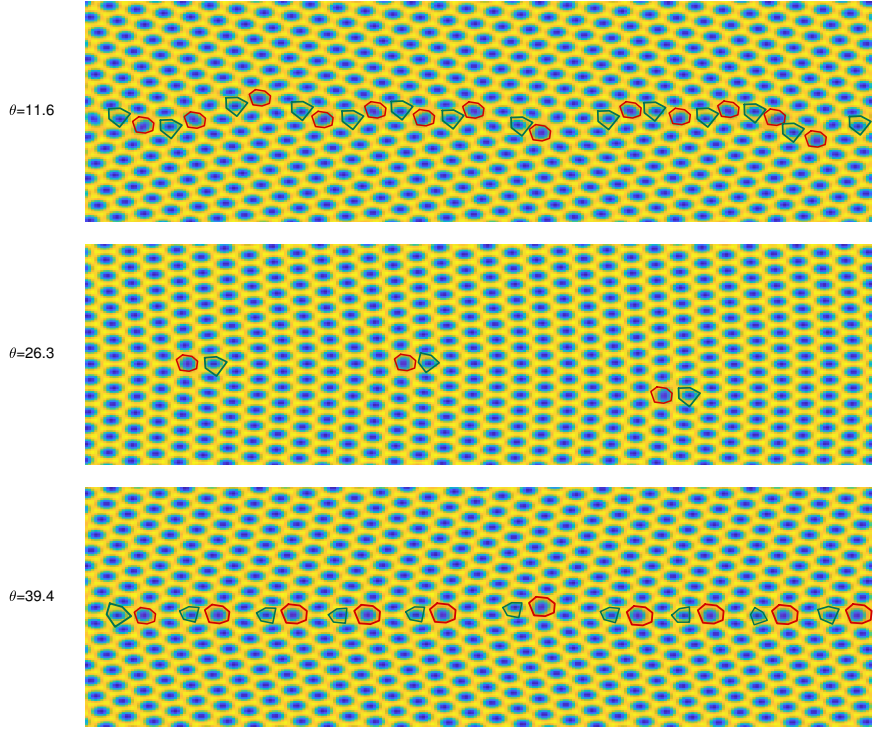

Fig. 5. (Example 4.3) Snapshots of the grain boundary mismatch of $\theta = 11.6°$(top), $\theta = 26.3°$(middle) and $\theta = 39.4°$(bottom). The green pentagons and red heptagons are located at the 5|7 dislocation dipoles.

We conduct two experiments in this example. The first one starts from a hexagonal lattice as the initial value with $\epsilon = 0.49$, while the other one starts from a square lattice as the initial value with $\epsilon = 0.017$. In both experiments, we solve NPFC models by proposed SAV scheme with discretizations of Fourier spectral method in space and BDF2 in time.

The phase evolutions are shown in Figs.6-7. We observe that the interface phase among bulk phase becomes from diffusive ones initially to sharp ones finally. These are exactly the discontinuities we have expected in this case. We believe that these sharp interface phenomenon have important physical meaning, since we only choose suitable kernels without any other extra conditions to achieve the sharp interface. Meanwhile, with the help of the hybrid algorithm provided in [21], the computation complexities of NPFC models are almost the same as LPFC models.

For concrete proof of such a discontinuity with integrable kernels, we present the following theorem for one-dimensional case. In two-dimensional case, one may observe the intersection of the graph of numerical solutions and plane $Y = 0$, in which it reduces to one-dimensional case.

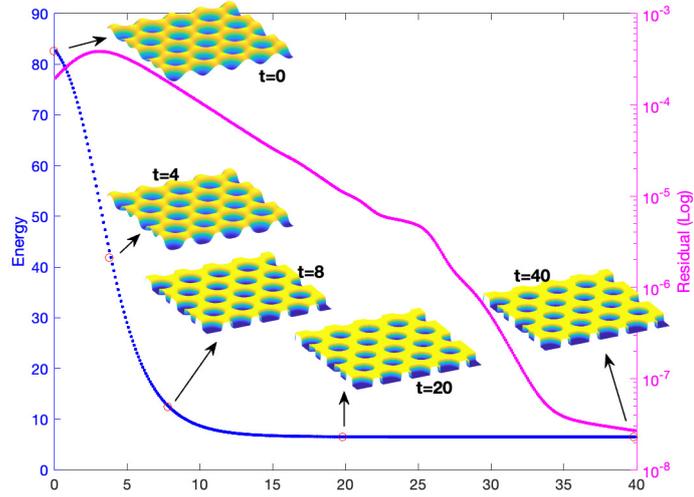

Fig. 6. (Example 4.4) Phase evolution for NPFC models starting from initial hexagonal lattice.

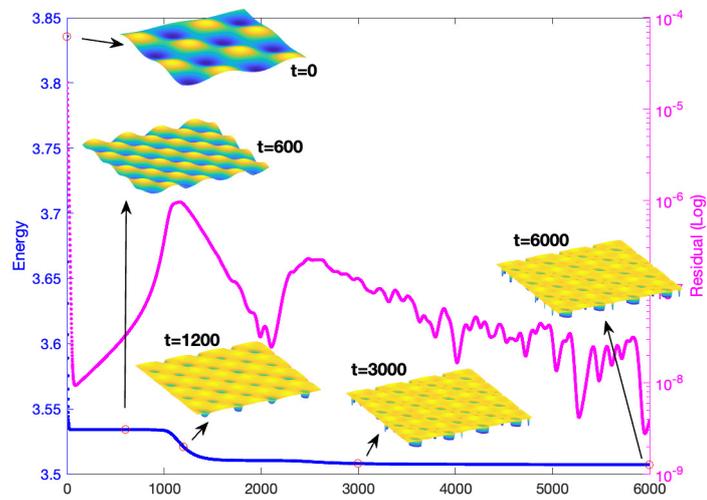

Fig. 7. (Example 4.4) Phase evolution for NPFC models starting from initial square lattice.

**Theorem 4.3.** *Assume $\phi^*$ is a solution to the aforementioned conserved NPFC model* (4.4) *and $\phi^*(x) = 0$ at point $x = x^*$. Then $\phi^*$ is continuous at $x = x^*$ when $|c_\delta - 1| \geq \sqrt{\epsilon}$; $\phi^*$ is discontinuous at $x = x^*$ when $|c_\delta - 1| < \sqrt{\epsilon}$ and $\rho_\delta + 2 - \frac{2}{c_\delta} < 0$.*

**Proof.** Without loss of generality we set $c_0 = 0$. Eq.(4.6) can be rewritten as

$$\phi^*(x)\left((\phi^*(x))^2 + (c_\delta - 1)^2 - \epsilon\right)$$
$$= -[\rho_\delta \star (\rho_\delta \star \phi^*)](x) + 2c_\delta[\rho_\delta \star \phi^*](x) - 2[\rho_\delta \star \phi^*](x)$$
$$:= I(x).$$

Note that

$$I(x) = -\left[\left(\rho_\delta + 2 - \frac{2}{c_\delta}\right) \star \rho_\delta \star \phi^*\right](x).$$

Assume $\phi^*$ is increasing over $[x^* - 2\delta, x^* + 2\delta]$, since $\phi^*(x^*) = 0$, then $\phi^*(x) \geq 0$ for any $x \in (x^*, x^* + \delta)$. Therefore $I(x^*) = 0$ and

$$I(x) = I(x) - I(x^*) > 0$$

if $\rho_\delta + 2 - 2/c_\delta < 0$. With this setting, when $|c_\delta - 1| \geq \sqrt{\epsilon}$, $\phi^*$ is continuous at $x = x^*$; when $|c_\delta - 1| < \sqrt{\epsilon}$, it requires that $(\phi^*(x))^2 + (c_\delta - 1)^2 - \epsilon \geq 0$, that is

$$\phi^*(x) \geq \sqrt{\epsilon - (c_\delta - 1)^2} > 0,$$

which yields discontinuity of $\phi^*(x)$ at $x = x^*$ since $\phi^*(x^*) = 0$. □

**Remark 4.3.** Restriction $\rho_\delta + 2 - 2/c_\delta < 0$ is not necessary in numerical examples. Its introduction is just for the ease of proof. The existence of discontinuous solutions established here is similar to that given in [20]. We also note a related study of discontinuous solutions to the nonlocal Cahn-Hilliard equation in [13].

The following example is designated to verify Theorem 4.3.

**Example 4.5.** We verify Theorem 4.3 over the domain $[0, 10]$ by choosing constant integrable kernel $\rho_\delta = \frac{3}{\delta^3}$ in one-dimensional case. Parameter $\epsilon = 1/4$. The initial value is given by $\phi_0 = \sin\left(\frac{2}{5}\pi x\right)$. We plot numerical solutions in Fig. 8 for different values of $\delta$ at $T = 500$.

Note that in this case $c_\delta = \frac{6}{\delta^2}$. We see from Fig.8 that with $\delta = 1.5$ in which case $\delta$ satisfies the condition $|c_\delta - 1| > \sqrt{\epsilon}$, the graph of numerical solution is continuous; with $\delta$ increasing to 2.5, discontinuity appears between phases. In this case $\delta$ satisfies $|c_\delta - 1| < \sqrt{\epsilon}$ while $\rho_\delta + 2 - 2/c_\delta > 0$. For $\delta = 2.6$ which satisfies $|c_\delta - 1| < \sqrt{\epsilon}$ and $\rho_\delta + 2 - 2/c_\delta < 0$, the graph is still discontinuous. With $\delta = 4.5$ satisfying $|c_\delta - 1| > \sqrt{\epsilon}$, we see the graph of numerical solution deforms to continuous interface again. This example confirms the findings of Theorem 4.3.

## 5. Concluding remarks

The present work presents a new data-driven 2D NPFC models. We first propose a nonlocal analogue free energy with the ND operators, which are more general than the original local PFC and also can encompass other fractional variants. Our studies show that the NPFC can provide much better fitting of the structure factor obtained by experiment data, thus enhancing its modeling capability.

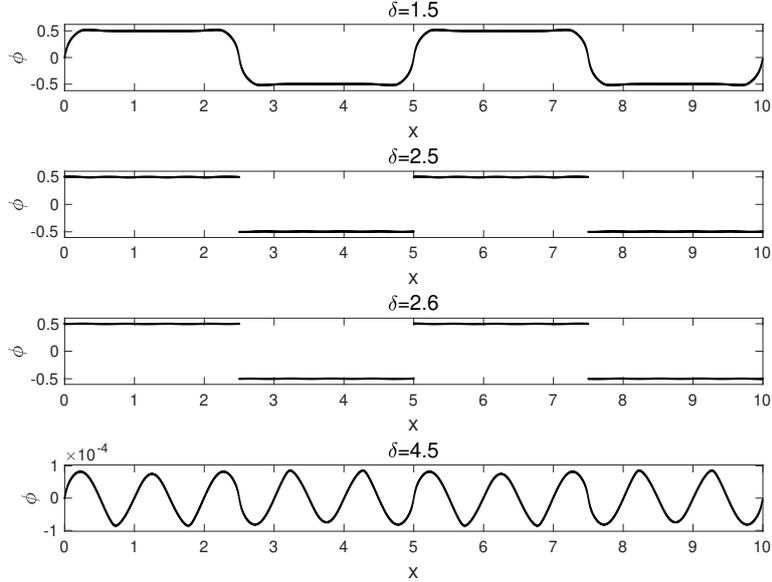

Fig. 8. (Example 4.5) Graph of numerical solutions for various values taken by $\delta$.

The NPFC model can be viewed as the gradient flow in the dual space of the energy space associated with the ND operator. We numerically solve NPFC models using an SAV scheme with Fourier spectral methods in space and BDF2 in time, together with the hybrid algorithm proposed in [21] to evaluate Fourier symbols of ND operator with specially chosen kernels. Numerical examples show that the proposed NPFC models, in particular those with sign-changing kernels that incorporate both repulsive and attractive interactions, can offer significant advantages in modeling the liquid-solid transition systems, such as in characterizing elastic and plastic deformations, as well as the anisotropy of solid-solid and solid-liquid interfaces.

While the current work is focused on the 2D case, the model can in principle be extended to three dimensions, which could be of greater relevance to applications in materials science. More careful studies in this direction will be pursued in the future.

## Acknowledgment

The work of Q. Du is supported in part by the NSF DMS 2309245 and DMS-1937254. The work of J. Yang is supported by the National Science Foundation of China(NSFC-12271240), and the Shenzhen Natural Science Fund (RCJC20210609103819018).

**Appendix**

### 5.1. *Proof of Lemma 3.1*

**Proof.** We have, from the energy stability, that

$$\|\phi_N^n\| + \|\mathcal{L}_\delta \phi_N^n\| + |r^n| \leq C \quad \text{since} \quad \|\mathcal{L}_\delta \phi_0\| \leq \|\Delta \phi_0\| \leq C.$$

By substituting the second equation in Eq.(3.10) into the first one, we obtain

$$\left(\frac{3\phi_N^{n+1} - 4\phi_N^n + \phi_N^{n-1}}{2\tau}, q\right) = (\tilde{\mathcal{L}}_\delta \mathcal{G}_\delta \phi_N^{n+1}, q) + \beta(\tilde{\mathcal{L}}_\delta \phi_N^{n+1}, q) + r^{n+1}\left(\tilde{\mathcal{L}}_\delta H(\bar{\phi}_N^{n+1}), q\right).$$

Testing the above equation with $q = \Delta^2 \phi_N^{n+1}$ to yield

$$\frac{1}{4\tau}\left(\|\Delta \phi_N^{n+1}\|^2 - \|\Delta \phi_N^n\|^2 + \|\Delta(2\phi_N^{n+1} - \phi_N^n)\|^2\right.$$
$$\left. - \|\Delta(2\phi_N^n - \phi_N^{n-1})\|^2 + \|\Delta(\phi_N^{n+1} - 2\phi_N^n + \phi_N^{n-1})\|^2\right)$$
$$= (\tilde{\mathcal{L}}_\delta \mathcal{G}_\delta \Delta \phi_N^{n+1}, \Delta \phi_N^{n+1}) + \beta(\tilde{\mathcal{L}}_\delta \Delta \phi_N^{n+1}, \Delta \phi_N^{n+1}) + r^{n+1}\left(\Delta H(\bar{\phi}_N^{n+1}), \tilde{\mathcal{L}}_\delta \Delta \phi_N^{n+1}\right)$$
$$\leq (\tilde{\mathcal{L}}_\delta \mathcal{G}_\delta \Delta \phi_N^{n+1}, \Delta \phi_N^{n+1}) + \beta(\tilde{\mathcal{L}}_\delta \Delta \phi_N^{n+1}, \Delta \phi_N^{n+1}) + C\left(\Delta H(\bar{\phi}_N^{n+1}), \tilde{\mathcal{L}}_\delta \Delta \phi_N^{n+1}\right).$$

With $\mathcal{G}_\delta = (\mathcal{L}_\delta + I)^2$, the above inequality can be further evaluated as

$$\frac{1}{4\tau}\left(\|\Delta \phi_N^{n+1}\|^2 - \|\Delta \phi_N^n\|^2 + \|\Delta(2\phi_N^{n+1} - \phi_N^n)\|^2\right.$$
$$\left. - \|\Delta(2\phi_N^n - \phi_N^{n-1})\|^2 + \|\Delta(\phi_N^{n+1} - 2\phi_N^n + \phi_N^{n-1})\|^2\right)$$
$$= \left(\mathcal{L}_\delta\left(\frac{3}{4}\mathcal{L}_\delta^2 + \left(\frac{1}{2}\mathcal{L}_\delta + 2\right)^2\right)\Delta \phi_N^{n+1}, \Delta \phi_N^{n+1}\right)$$
$$+ ((\beta - 3)\mathcal{L}_\delta \Delta \phi_N^{n+1}, \Delta \phi_N^{n+1}) + C\left(\Delta H(\bar{\phi}_N^{n+1}), \tilde{\mathcal{L}}_\delta \Delta \phi_N^{n+1}\right)$$
$$\leq \left(\frac{3}{4}\mathcal{L}_\delta^3 \Delta \phi_N^{n+1}, \Delta \phi_N^{n+1}\right) + C\left(\frac{1}{2}\|\Delta H(\bar{\phi}_N^{n+1})\|^2 + \frac{1}{2}\|\mathcal{L}_\delta \Delta \phi_N^{n+1}\|^2\right) := I_r.$$

where $\beta > 3$ is chosen and the non-negative definite property of the operator $\mathcal{L}_\delta$ is used. By $\tilde{\mathcal{L}}_\delta = \mathcal{L}_\delta$, $\Delta U(\phi) = (U''(\phi)|\nabla \phi|^2 + U'(\phi)\Delta \phi)$, equation (3.9) and interpolation inequality, it follows that

$$\|\Delta H(\phi)\| \leq C(\|U''\|_\infty \|\nabla \phi\|_{L^4}^2 + \|U'\|_\infty \|\Delta \phi\|)$$
$$\leq C\|\nabla \phi\| \|\Delta \phi\| + C\|\Delta \phi\| \leq C\|\Delta \phi\|,$$

in which the assumption (3.4) ensures the $L^\infty$-bound of $U'$ and $U''$. We then then simplify $I_r$ into

$$I_r \leq \left(\left(\frac{3}{4}\mathcal{L}_\delta^3 + \frac{C}{2}\right)\Delta \phi_N^{n+1}, \Delta \phi_N^{n+1}\right) + \frac{C}{2}\left(\|\Delta \phi_N^n\|^2 + \|\Delta \phi_N^{n-1}\|^2\right).$$

With the constant $C$ as appeared in the above, for the function

$$g(\lambda_\delta) = \frac{3}{4}\lambda_\delta^3 + \frac{C}{2}\lambda_\delta^2 = \frac{3}{4}\lambda_\delta^2\left(\lambda_\delta + \frac{2}{3}C\right),$$

we can note that $g(\lambda_\delta) \leq 0$ for $\lambda_\delta \leq -\frac{2}{3}C$. And for $-\frac{2}{3}C < \lambda_\delta < 0$, $g(\lambda_\delta)$ is bounded since $g$ is continuous and $g\left(-\frac{2}{3}C\right) = 0, g(0) = 0$. Thus

$$I_r \leq C(\|\Delta\phi_N^{n+1}\|^2 + \|\Delta\phi_N^n\|^2 + \|\Delta\phi_N^{n-1}\|^2).$$

Combining the above yields

$$\frac{1}{4\tau}\left(\|\Delta\phi_N^{n+1}\|^2 - \|\Delta\phi_N^n\|^2 + \|\Delta(2\phi_N^{n+1} - \phi_N^n)\|^2\right.$$
$$\left. - \|\Delta(2\phi_N^n - \phi_N^{n-1})\|^2 + \|\Delta(\phi_N^{n+1} - 2\phi_N^n + \phi_N^{n-1})\|^2\right)$$
$$\leq C(\|\Delta\phi_N^{n+1}\|^2 + \|\Delta\phi_N^n\|^2 + \|\Delta\phi_N^{n-1}\|^2),$$

from which $\|\Delta\phi_N^n\| \leq C$ follows with the help of Grownwall's inequality. Due to Sobolev embedding theorem $H^{1+\epsilon} \subset L^\infty$ for any $\epsilon > 0$ in two dimensions, we have the desired estimate $\|\phi_N^n\|_{L^\infty} \leq \|\phi_N^n\|_{H^2} \leq C$. □

### 5.2. *Proof of Lemma 3.2*

**Proof.** Without loss of generality, let us assume the Lipschitz conntinuity of $v$. Denote by

$$L(x, y) = \frac{g(x) - g(y)}{x - y}.$$

Since $g$ is gloablly Lipschitz continuous from the assumption, $L(x, y)$ is thus bounded. Note the following splitting

$$g(u(y)) - g(v(y)) - (g(u(x)) - g(v(x)))$$
$$= g(u(y)) - g(u(x)) - (g(v(y)) - g(v(x)))$$
$$= (u(y) - u(x))L(u(y), u(x)) - (v(y) - v(x))L(v(y), v(x))$$
$$= (u(y) - u(x))L(u(y), u(x)) - (v(y) - v(x))L(u(y), u(x))$$
$$\quad + (v(y) - v(x))L(u(y), u(x)) - (v(y) - v(x))L(v(y), v(x))$$
$$= \{(e(y) - e(x))L(u(y), u(x))\} + \{(v(y) - v(x))(L(u(y), u(x)) - L(v(y), u(x)))\}$$
$$\quad + \{(v(y) - v(x))(L(v(y), u(x)) - L(v(y), v(x)))\}$$
$$:= I_1 + I_2 + I_3.$$

With the first term $I_1 = (e(y) - e(x))L(u(y), u(x))$, it holds that

$$\iint_\Omega \rho_\delta(|x - y|)|I_1|^2 dy dx \leq C \iint_\Omega \rho_\delta(|x - y|)(e(y) - e(x))^2 dy dx = C(-\mathcal{L}_\delta e, e).$$

To estimate the last two terms $I_2$ and $I_3$, we first note that

$$L'_u(u, v) := \frac{\partial L(u, v)}{\partial u} = \frac{g'(u)(u - v) - (g(u) - g(v))}{(u - v)^2} = \frac{g'(u) - g'(\zeta)}{u - v}$$

for $\zeta$ in between $u$ and $v$. By the Lipschitz continuity of $g'$, we thus see that

$$|L'_u(u, v)| \leq C\frac{|u - \zeta|}{|u - v|} \leq C,$$

which implies $L(u,v)$ is globally Lipschitz continuous in $u$. Similarly, $L(u,v)$ is also globally Lipschitz continuous in $v$. Hence, the last two terms can be estimated as

$$\iint_\Omega \rho_\delta(|x-y|)|I_2|^2 dy dx$$
$$= \iint_\Omega \rho_\delta(|x-y|)|(v(y)-v(x))(L(u(y),u(x))-L(v(y),u(x)))|^2 dy dx$$
$$\leq C \iint_\Omega \rho_\delta(|x-y|)(v(y)-v(x))^2 |e(y)|^2 dy dx$$
$$= C \int_\Omega e^2(y) \int_\Omega \rho_\delta(|y-x|)|y-x|^2 dx dy$$
$$= C(e,e)$$

and

$$\iint_\Omega \rho_\delta(|x-y|)|I_3|^2 dy dx$$
$$= \iint_\Omega \rho_\delta(|x-y|)|(v(y)-v(x))(L(v(y),u(x))-L(v(y),v(x)))|^2 dy dx$$
$$\leq C \iint_\Omega \rho_\delta(|x-y|)(v(y)-v(x))^2 |e(x)|^2 dy dx$$
$$= C \int_\Omega e^2(x) \int_\Omega \rho_\delta(|y-x|)|y-x|^2 dy dx$$
$$= C(e,e)$$

where the second moment condition $\int_{\mathcal{B}_\delta(0)} \rho_\delta(|s|)s^2 ds = C$ and Lipschitz continuity of $v$ are used. Finally the desired estimate is derived by

$$(-\mathcal{L}_\delta(g(u)-g(v)), g(u)-g(v))$$
$$= \frac{1}{2} \iint_\Omega \rho_\delta(|x-y|)(g(u(y))-g(v(y))-(g(u(x))-g(v(x))))^2 dy dx$$
$$\leq C \sum_{k=1}^3 \iint_\Omega \rho_\delta(|x-y|)|I_k|^2 dy dx$$
$$\leq C((-\mathcal{L}_\delta e, e)+(e,e)).$$

From the symmetry, we see that the same proof works with $u$ being Lipschitz if a different splitting is used. In turn, we can also do the splitting with a convex linear combination of $u$ and $v$. □

### 5.3. *Proof of Theorem 3.2*

**Proof.** For simplicity, we set

$$e_\phi^n = \phi_N^n - \Pi_h \phi(t_n) + \Pi_h \phi(t_n) - \phi(t_n) = \bar{e}_\phi^n + \check{e}_\phi^n,$$
$$e_\mu^n = \mu_N^n - \Pi_h \mu(t_n) + \Pi_h \mu(t_n) - \mu(t_n) = \bar{e}_\mu^n + \check{e}_\mu^n,$$

$$e_r^n = r^n - r(t_n)$$

where $\Pi_h$ is spectral orthogonal projection defined by Eq.(3.7). Subtracting Eq.(3.8) from Eq.(3.10) at $t = t_{n+1}$, we have

$$\begin{cases} (3\bar{e}_\phi^{n+1} - 4\bar{e}_\phi^n + \bar{e}_\phi^{n-1}, q) = 2\tau(\tilde{\mathcal{L}}_\delta \bar{e}_\mu^{n+1}, q) + (Q_1^{n+1}, q) & \forall q \in S_N, \\ (\bar{e}_\mu^{n+1}, w) = (\mathcal{G}_\delta \bar{e}_\phi^{n+1}, w) + (\beta \bar{e}_\phi^{n+1}, w) \\ \qquad + (r^{n+1} H(\bar{\phi}_N^{n+1}), w) - (r(t_{n+1}) H(\phi(t_{n+1})), w) & \forall w \in S_N, \\ 3e_r^{n+1} - 4e_r^n + e_r^{n-1} = \dfrac{1}{2}(H(\bar{\phi}_N^{n+1}), 3\bar{e}_\phi^{n+1} - 4\bar{e}_\phi^n + \bar{e}_\phi^{n-1}) + Q_2^{n+1} - Q_3^{n+1} \\ \qquad + \dfrac{1}{2}(H(\bar{\phi}_N^{n+1}) - H(\phi(t_{n+1})), 3\phi(t_{n+1}) - 4\phi(t_n) + \phi(t_{n-1})) \end{cases}$$

where $\bar{\phi}_N^{n+1}$ is the extrapolation as defined in Eq.(3.11) and

$$\begin{aligned} Q_1^{n+1} &= 2\tau \partial_t \phi(t_{n+1}) - (3\phi(t_{n+1}) - 4\phi(t_n) + \phi(t_{n-1})) \\ &= 2\int_{t_{n+1}}^{t_n} (t_n - s)^2 \partial_t^3 \phi(s) ds - \frac{1}{2}\int_{t_{n+1}}^{t_{n-1}} (t_{n-1} - s)^2 \partial_t^3 \phi(s) ds, \\ Q_2^{n+1} &= 2\tau d_t r(t_{n+1}) - (3r(t_{n+1}) - 4r(t_n) + r(t_{n-1})) \\ &= 2\int_{t_{n+1}}^{t_n} (t_n - s)^2 d_t^3 r(s) ds - \frac{1}{2}\int_{t_{n+1}}^{t_{n-1}} (t_{n-1} - s)^2 d_t^3 r(s) ds, \\ Q_3^{n+1} &= \frac{1}{2}(H(\phi(t_{n+1})), Q_1^{n+1}). \end{aligned}$$

Taking $q = \bar{e}_\mu^{n+1}$, $w = 3\bar{e}_\phi^{n+1} - 4\bar{e}_\phi^n + \bar{e}_\phi^{n-1}$ as the test functions in the above and taking the inner product of the equation for $e_r^{n+1}$ with $2e_r^{n+1}$, we get

$$\begin{aligned} (3\bar{e}_\phi^{n+1} - 4\bar{e}_\phi^n + \bar{e}_\phi^{n-1}, \bar{e}_\mu^{n+1}) + 2\tau \|(-\tilde{\mathcal{L}}_\delta)^{\frac{1}{2}} \bar{e}_\mu^{n+1}\|^2 &= (Q_1^{n+1}, \bar{e}_\mu^{n+1}), & (5.1) \\ (\bar{e}_\mu^{n+1}, 3\bar{e}_\phi^{n+1} - 4\bar{e}_\phi^n + \bar{e}_\phi^{n-1}) &= (\mathcal{G}_\delta \bar{e}_\phi^{n+1}, 3\bar{e}_\phi^{n+1} - 4\bar{e}_\phi^n + \bar{e}_\phi^{n-1}) \\ &\quad + \beta(\bar{e}_\phi^{n+1}, 3\bar{e}_\phi^{n+1} - 4\bar{e}_\phi^n + \bar{e}_\phi^{n-1}) + (r^{n+1} H(\bar{\phi}_N^{n+1}), 3\bar{e}_\phi^{n+1} - 4\bar{e}_\phi^n + \bar{e}_\phi^{n-1}) \\ &\quad - (r(t_{n+1}) H(\phi(t_{n+1})), 3\bar{e}_\phi^{n+1} - 4\bar{e}_\phi^n + \bar{e}_\phi^{n-1}), & (5.2) \\ 2e_r^{n+1}(3e_r^{n+1} - 4e_r^n + e_r^{n-1}) &= e_r^{n+1}(H(\bar{\phi}_N^{n+1}), 3\bar{e}_\phi^{n+1} - 4\bar{e}_\phi^n + \bar{e}_\phi^{n-1}) \\ &\quad + 2e_r^{n+1}(Q_2^{n+1} - Q_3^{n+1}) \\ &\quad + e_r^{n+1}(H(\bar{\phi}_N^{n+1}) - H(\phi(t_{n+1})), 3\phi(t_{n+1}) - 4\phi(t_n) + \phi(t_{n-1})). & (5.3) \end{aligned}$$

The term on the right-hand side of Eq.(5.1) can be estimated by

$$\begin{aligned} (Q_1^{n+1}, \bar{e}_\mu^{n+1}) &= ((-\tilde{\mathcal{L}}_\delta)^{-\frac{1}{2}} Q_1^{n+1}, (-\tilde{\mathcal{L}}_\delta)^{\frac{1}{2}} \bar{e}_\mu^{n+1}) \\ &\leq \frac{\tau}{4} \|(-\tilde{\mathcal{L}}_\delta)^{\frac{1}{2}} \bar{e}_\mu^{n+1}\|^2 + \frac{1}{\tau} \|(-\tilde{\mathcal{L}}_\delta)^{-\frac{1}{2}} Q_1^{n+1}\|^2 \\ &\leq \frac{\tau}{4} \|(-\tilde{\mathcal{L}}_\delta)^{\frac{1}{2}} \bar{e}_\mu^{n+1}\|^2 + C\tau^4 \int_{t_{n-1}}^{t_{n+1}} \|(-\tilde{\mathcal{L}}_\delta)^{-\frac{1}{2}} \partial_t^3 \phi(s)\|^2 ds. \end{aligned}$$

On the right-hand side of Eq.(5.2), the first two terms can be estimated using Eq.(3.12), i.e.,

$$(\mathcal{G}_\delta \bar{e}_\phi^{n+1}, 3\bar{e}_\phi^{n+1} - 4\bar{e}_\phi^n + \bar{e}_\phi^{n-1}) + \beta(\bar{e}_\phi^{n+1}, 3\bar{e}_\phi^{n+1} - 4\bar{e}_\phi^n + \bar{e}_\phi^{n-1})$$
$$= \frac{\beta}{2}(\|\bar{e}_\phi^{n+1}\|^2 - \|\bar{e}_\phi^n\|^2 + \|2\bar{e}_\phi^{n+1} - \bar{e}_\phi^n\|^2 + \|\bar{e}_\phi^{n+1} - 2\bar{e}_\phi^n + \bar{e}_\phi^{n-1}\|^2 - \|2\bar{e}_\phi^n - \bar{e}_\phi^{n-1}\|^2)$$
$$+ \frac{1}{2}\Big(\|(I+\mathcal{L}_\delta)\bar{e}_\phi^{n+1}\|^2 - \|(I+\mathcal{L}_\delta)\bar{e}_\phi^n\|^2 + \|(I+\mathcal{L}_\delta)(2\bar{e}_\phi^{n+1} - \bar{e}_\phi^n)\|^2$$
$$+ \|(I+\mathcal{L}_\delta)(\bar{e}_\phi^{n+1} - 2\bar{e}_\phi^n + \bar{e}_\phi^{n-1})\|^2 - \|(I+\mathcal{L}_\delta)(2\bar{e}_\phi^n - \bar{e}_\phi^{n-1})\|^2\Big)$$

and the last two terms can be transformed to

$$(r^{n+1}H(\bar{\phi}_N^{n+1}), 3\bar{e}_\phi^{n+1} - 4\bar{e}_\phi^n + \bar{e}_\phi^{n-1}) - (r(t_{n+1})H(\phi(t_{n+1})), 3\bar{e}_\phi^{n+1} - 4\bar{e}_\phi^n + \bar{e}_\phi^{n-1})$$
$$= (e_r^{n+1}H(\bar{\phi}_N^{n+1}), 3\bar{e}_\phi^{n+1} - 4\bar{e}_\phi^n + \bar{e}_\phi^{n-1})$$
$$+ r(t_{n+1})(H(\bar{\phi}_N^{n+1}) - H(\phi(t_{n+1})), 3\bar{e}_\phi^{n+1} - 4\bar{e}_\phi^n + \bar{e}_\phi^{n-1}).$$

Due to the boundedness of numerical solution (3.13), we have $\|\bar{e}_\phi^n\|_{L^\infty}$ is bounded and thus the last term on the right-hand side of the above equation can be controlled as

$$r(t_{n+1})(H(\bar{\phi}_N^{n+1}) - H(\phi(t_{n+1})), 3\bar{e}_\phi^{n+1} - 4\bar{e}_\phi^n + \bar{e}_\phi^{n-1})$$
$$= r(t_{n+1})\left(H(\bar{\phi}_N^{n+1}) - H(\phi(t_{n+1})), 2\tau\tilde{\mathcal{L}}_\delta \bar{e}_\mu^{n+1}\right)$$
$$+ r(t_{n+1})\left(H(\bar{\phi}_N^{n+1}) - H(\phi(t_{n+1})), Q_1^{n+1}\right)$$
$$\leq \frac{\tau}{4}\|(-\tilde{\mathcal{L}}_\delta)^{\frac{1}{2}}\bar{e}_\mu^{n+1}\|^2 + C\tau\|(-\tilde{\mathcal{L}}_\delta)^{\frac{1}{2}}H(\bar{\phi}_N^{n+1}) - (-\tilde{\mathcal{L}}_\delta)^{\frac{1}{2}}H(\phi(t_{n+1}))\|^2$$
$$+ \frac{C}{\tau}\|(-\tilde{\mathcal{L}}_\delta)^{-\frac{1}{2}}Q_1^{n+1}\|^2.$$

Note that with

$$\mathcal{W}(\phi) = \int_\Omega \{F(\phi)\} - \frac{\beta}{2}\phi^2 + C_H,$$

then

$$H(\bar{\phi}_N^{n+1}) - H(\phi(t_{n+1})) = \frac{U(\bar{\phi}_N^{n+1})}{\sqrt{\mathcal{W}(\bar{\phi}_N^{n+1})}} - \frac{U(\phi(t_{n+1}))}{\sqrt{\mathcal{W}(\phi(t_{n+1}))}}$$
$$= \frac{U(\bar{\phi}_N^{n+1}) - U(\phi(t_{n+1}))}{\sqrt{\mathcal{W}(\bar{\phi}_N^{n+1})}} + \frac{U(\phi(t_{n+1}))}{\sqrt{\mathcal{W}(\bar{\phi}_N^{n+1})}\sqrt{\mathcal{W}(\phi(t_{n+1}))}} \frac{\mathcal{W}(\bar{\phi}_N^{n+1}) - \mathcal{W}(\phi(t_{n+1}))}{\sqrt{\mathcal{W}(\bar{\phi}_N^{n+1})} + \sqrt{\mathcal{W}(\phi(t_{n+1}))}}.$$

Since $\phi_N^n$ is bounded due to Lemma 3.1 with $\tilde{\mathcal{L}}_\delta = \mathcal{L}_\delta$, then $\left|\sqrt{\mathcal{W}(\bar{\phi}_N^{n+1})}\right|$, $\left|\sqrt{\mathcal{W}(\phi(t_{n+1}))}\right|$, $|U(\phi(t_{n+1}))|$ are bounded, and

$$\bar{\phi}_N^{n+1} - \phi(t_{n+1}) = 2e_\phi^n - e_\phi^{n-1} + \int_{t_{n-1}}^{t_{n+1}}(t-s)\partial_t^2\phi(s)ds,$$

using Lemma 3.2 with above splitting yields

$$\|(-\mathcal{L}_\delta)^{\frac{1}{2}}(H(\bar{\phi}_N^{n+1}) - H(\phi(t_{n+1})))\|^2$$
$$=(-\mathcal{L}_\delta(H(\bar{\phi}_N^{n+1}) - H(\phi(t_{n+1}))), H(\bar{\phi}_N^{n+1}) - H(\phi(t_{n+1})))$$
$$\leq C(-\mathcal{L}_\delta(U(\bar{\phi}_N^{n+1}) - U(\phi(t_{n+1}))), U(\bar{\phi}_N^{n+1}) - U(\phi(t_{n+1})))$$
$$+ C(-\mathcal{L}_\delta(\mathcal{W}(\bar{\phi}_N^{n+1}) - \mathcal{W}(\phi(t_{n+1}))), \mathcal{W}(\bar{\phi}_N^{n+1}) - \mathcal{W}(\phi(t_{n+1})))$$
$$\leq C(-\mathcal{L}_\delta(2e_\phi^n - e_\phi^{n-1}), 2e_\phi^n - e_\phi^{n-1}) + C(2e_\phi^n - e_\phi^{n-1}, 2e_\phi^n - e_\phi^{n-1})$$
$$+ C\left(-\mathcal{L}_\delta \int_{t_{n-1}}^{t_{n+1}} (t-s)\partial_t^2\phi(s)ds, \int_{t_{n-1}}^{t_{n+1}} (t-s)\partial_t^2\phi(s)ds\right)$$
$$+ C\left(\int_{t_{n-1}}^{t_{n+1}} (t-s)\partial_t^2\phi(s)ds, \int_{t_{n-1}}^{t_{n+1}} (t-s)\partial_t^2\phi(s)ds\right)$$
$$\leq C(\|(I+\mathcal{L}_\delta)(2e_\phi^n - e_\phi^{n-1}))\|^2 + \|2e_\phi^n - e_\phi^{n-1}\|^2)$$
$$+ C\tau^4 \int_{t_{n-1}}^{t_{n+1}} (\|(-\mathcal{L}_\delta)^{\frac{1}{2}}\partial_t^2\phi(s)\|^2 + \|\partial_t^2\phi(s)\|^2)ds$$

where the last inequality follows from the the following inequality: for any $u$,

$$\|(-\mathcal{L}_\delta)^{\frac{1}{2}}u\|^2 = ((-\mathcal{L}_\delta)u, u) = (-(\mathcal{L}_\delta + I)u, u) + (u, u)$$
$$\leq C\left(\|(I+\mathcal{L}_\delta)u\|^2 + \|u\|^2\right). \tag{5.4}$$

Thus

$$r(t_{n+1})(H(\bar{\phi}_N^{n+1}) - H(\phi(t_{n+1})), 3\bar{e}_\phi^{n+1} - 4\bar{e}_\phi^n + \bar{e}_\phi^{n-1})$$
$$\leq \frac{\tau}{4}\|(-\mathcal{L}_\delta)^{\frac{1}{2}}\bar{e}_\mu^{n+1}\|^2 + C\tau(\|2\bar{e}_\phi^n - \bar{e}_\phi^{n-1}\|^2 + \|(-\mathcal{L}_\delta)^{\frac{1}{2}}(2\bar{e}_\phi^n - \bar{e}_\phi^{n-1})\|^2)$$
$$+ C\tau^4 \int_{t_{n-1}}^{t_{n+1}} (\|(-\mathcal{L}_\delta)^{\frac{1}{2}}\partial_t^2\phi(s)\|^2 + \|\partial_t^2\phi(s)\|^2)ds$$
$$+ C\tau^4 \int_{t_{n-1}}^{t_{n+1}} \|(-\mathcal{L}_\delta)^{-\frac{1}{2}}\partial_t^3\phi(s)\|^2 ds$$

where in the last inequality we use the Lipschitz continuity assumption on $\phi$ and Lemma 3.2 which holds due to the assumption (3.4) and equation (3.9). $Q_2^{n+1}$ is, similar to $Q_1^{n+1}$, estimated as

$$|Q_2^{n+1}|^2 \leq C\tau^5 \int_{t_{n-1}}^{t_{n+1}} |\mathrm{d}_t^3 r(s)|^2 ds$$

and

$$|Q_3^{n+1}|^2 \leq C\tau^5 \|(-\mathcal{L}_\delta)^{\frac{1}{2}}H(\phi(t_n))\|^2 \int_{t_{n-1}}^{t_{n+1}} \|(-\mathcal{L}_\delta)^{-\frac{1}{2}}\partial_t^3\phi(s)\|^2 ds$$
$$\leq C\tau^5 \int_{t_{n-1}}^{t_{n+1}} \|(-\mathcal{L}_\delta)^{-\frac{1}{2}}\partial_t^3\phi(s)\|^2 ds.$$

Then the second term on the right-hand side of Eq.(5.3) can be estimated as

$$2e_r^{n+1}(Q_2^{n+1} - Q_3^{n+1}) \leq \tau |e_r^{n+1}|^2 + \frac{C}{\tau}(|Q_2^{n+1}|^2 + |Q_3^{n+1}|^2)$$
$$\leq \tau |e_r^{n+1}|^2 + C\tau^4 \int_{t_{n-1}}^{t_{n+1}} (|\mathrm{d}_t^3 r(s)|^2 + \|(-\mathcal{L}_\delta)^{-\frac{1}{2}} \partial_t^3 \phi(s)\|^2) ds.$$

The last term on the right-hand side of Eq.(5.3) can be controlled as

$$e_r^{n+1}(H(\bar{\phi}_N^{n+1}) - H(\phi(t_{n+1})), 3\phi(t_{n+1}) - 4\phi(t_n) + \phi(t_{n-1}))$$
$$= e_r^{n+1}((-\mathcal{L}_\delta)^{\frac{1}{2}}(H(\bar{\phi}_N^{n+1}) - H(\phi(t_{n+1}))), (-\mathcal{L}_\delta)^{-\frac{1}{2}}(3\phi(t_{n+1}) - 4\phi(t_n) + \phi(t_{n-1})))$$
$$= e_r^{n+1}\left((-\mathcal{L}_\delta)^{\frac{1}{2}}(H(\bar{\phi}_N^{n+1}) - H(\phi(t_{n+1}))), (-\mathcal{L}_\delta)^{-\frac{1}{2}}\left(\int_{t_n}^{t_{n+1}} \partial_t \phi(s) ds + \int_{t_n}^{t_{n-1}} \partial_t \phi(s) ds\right)\right)$$
$$\leq C\tau \|\partial_t \phi\|^2_{L^\infty((0,T];H((-\mathcal{L}_\delta)^{-\frac{1}{2}}))} (|e_r^{n+1}|^2 + \|(-\mathcal{L}_\delta)^{\frac{1}{2}} H(\bar{\phi}_N^{n+1}) - (-\mathcal{L}_\delta)^{\frac{1}{2}} H(\phi(t_{n+1}))\|^2)$$
$$\leq C\tau \|\partial_t \phi\|^2_{L^\infty((0,T];H((-\mathcal{L}_\delta)^{-\frac{1}{2}}))} (|e_r^{n+1}|^2 + \|2\bar{e}_\phi^n - \bar{e}_\phi^{n-1}\|^2 + \|(-\mathcal{L}_\delta)^{\frac{1}{2}}(2\bar{e}_\phi^n - \bar{e}_\phi^{n-1})\|^2)$$
$$+ C\tau^4 \|\partial_t \phi\|^2_{L^\infty((0,T];H((-\mathcal{L}_\delta)^{-\frac{1}{2}}))} \int_{t_{n-1}}^{t_{n+1}} (\|(-\mathcal{L}_\delta)^{\frac{1}{2}} \partial_t^2 \phi(s)\|^2 + \|\partial_t^2 \phi(s)\|^2) ds$$

where Lemma 3.2 is applied in the second inequality.

Combining the above equations we have

$$2\tau \|(-\mathcal{L}_\delta)^{\frac{1}{2}} \bar{e}_\mu^{n+1}\|^2 + |e_r^{n+1}|^2 - |e_r^n|^2$$
$$+ |2e_r^{n+1} - e_r^n|^2 + |e_r^{n+1} - 2e_r^n + e_r^{n-1}|^2 - |2e_r^n - e_r^{n-1}|^2$$
$$+ \frac{\beta}{2}(\|\bar{e}_\phi^{n+1}\|^2 - \|\bar{e}_\phi^n\|^2 + \|2\bar{e}_\phi^{n+1} - \bar{e}_\phi^n\|^2 + \|\bar{e}_\phi^{n+1} - 2\bar{e}_\phi^n + \bar{e}_\phi^{n-1}\|^2 - \|2\bar{e}_\phi^n - \bar{e}_\phi^{n-1}\|^2)$$
$$+ \frac{1}{2}\Big(\|(I + \mathcal{L}_\delta)\bar{e}_\phi^{n+1}\|^2 - \|(I + \mathcal{L}_\delta)\bar{e}_\phi^n\|^2 + \|(I + \mathcal{L}_\delta)(2\bar{e}_\phi^{n+1} - \bar{e}_\phi^n)\|^2$$
$$+ \|(I + \mathcal{L}_\delta)(\bar{e}_\phi^{n+1} - 2\bar{e}_\phi^n + \bar{e}_\phi^{n-1})\|^2 - \|(I + \mathcal{L}_\delta)(2\bar{e}_\phi^n - \bar{e}_\phi^{n-1})\|^2\Big)$$
$$\leq \frac{\tau}{2} \|(-\mathcal{L}_\delta)^{\frac{1}{2}} \bar{e}_\mu^{n+1}\|^2$$
$$+ C\tau^4 \int_{t_{n-1}}^{t_{n+1}} (\|(-\mathcal{L}_\delta)^{-\frac{1}{2}} \partial_t^3 \phi(s)\|^2 + \|(-\mathcal{L}_\delta)^{\frac{1}{2}} \partial_t^2 \phi(s)\|^2 + \|\partial_t^2 \phi(s)\|^2 + |\mathrm{d}_t^3 r(s)|^2) ds$$
$$+ C\tau (\|2\bar{e}_\phi^n - \bar{e}_\phi^{n-1}\|^2 + \|(-\mathcal{L}_\delta)^{\frac{1}{2}}(2\bar{e}_\phi^n - \bar{e}_\phi^{n-1})\|^2 + |e_r^{n+1}|^2),$$

which directly, by Gronwall's inequality and Eq.(5.4), results in

$$\|(I + \mathcal{L}_\delta)\bar{e}_\phi^{n+1}\|^2 + |e_r^{n+1}|^2 + \|\bar{e}_\phi^{n+1}\|^2$$
$$\leq C\tau^4 \int_0^T \left(\|(-\mathcal{L}_\delta)^{-\frac{1}{2}} \partial_t^3 \phi(s)\|^2 + \|(-\mathcal{L}_\delta)^{\frac{1}{2}} \partial_t^2 \phi(s)\|^2 + \|\partial_t^2 \phi(s)\|^2 + |\mathrm{d}_t^3 r(s)|^2\right) ds.$$

Finally, together with existing results on Fourier spectral method

$$\|\Pi_h u - u\|_k \leq C\|u\|_m h^{m-k}, \quad \forall 0 \leq k \leq m,$$

for $u \in H_{per}^m(\Omega)$ and

$$\|\Pi_h u - u\|_{L^\infty} \leq C e^{-c/h}$$

for analytic solution $u$, it holds that

$$\|\phi_N^n - \phi(t_n)\|^2 + |r^n - r(t_n)|^2 \leq C(h^{2m} + \tau^4)$$

for $\phi \in H_{per}^m(\Omega)$ and that

$$\|\phi_N^n - \phi(t_n)\|^2 + |r^n - r(t_n)|^2 \leq C(e^{-c/h} + \tau^4)$$

for analytic solution $\phi$, where constants $C$ and $c$ are independent of temporal step size and spatial step size. $\square$

### 5.4. *Proof of Lemma 4.1*

**Proof.** Since we have

$$\phi_{MN}^\delta - \phi_{MN}^0 = (\mathcal{L}_0^{-1} - \mathcal{L}_\delta^{-1})f_{MN},$$

then

$$\|\phi_{MN}^\delta - \phi_{MN}^0\|^2 = \sum_{\substack{|m|\leq M, |n|\leq N \\ m^2+n^2 \neq 0}} \left|\frac{1}{\lambda_\delta(m,n)} - \frac{1}{\lambda_0(m,n)}\right| |\hat{f}_{mn}|^2.$$

Thus it suffices to prove

$$\frac{1}{\delta^2}\left|\frac{1}{\lambda_\delta(m,n)} - \frac{1}{\lambda_0(m,n)}\right| := S_{mn} \leq C, \quad \forall |m| \leq M, |n| \leq N \text{ with } m^2 + n^2 \neq 0.$$

Due to $1 - \cos\theta \leq \frac{\theta^2}{2}$, we have

$$\begin{aligned}
0 \leq \delta^2 |\lambda_\delta(m,n)| &= \delta^2 \int_0^\delta r\rho_\delta(r) \int_0^{2\pi} \left(1 - \cos(r\sqrt{m^2+n^2}\cos\theta)\right) d\theta dr \\
&= \delta^2 \int_0^\delta r\delta^{-4}\rho(\frac{r}{\delta}) \int_0^{2\pi} \left(1 - \cos(r\sqrt{m^2+n^2}\cos\theta)\right) d\theta dr \\
&= \int_0^1 r\rho(r) \int_0^{2\pi} \left(1 - \cos(\delta r\sqrt{m^2+n^2}\cos\theta)\right) d\theta dr \\
&\leq \int_0^1 r\rho(r) \int_0^{2\pi} \left(\frac{\delta^2 r^2(m^2+n^2)\cos^2\theta}{2}\right) d\theta dr \\
&= \delta^2(m^2+n^2)\int_0^1 \rho(r)r^3 dr \int_0^{2\pi} \frac{\cos^2\theta}{2} d\theta \\
&= \delta^2(m^2+n^2),
\end{aligned}$$

where we used the form of the kernel and the second moment condition

$$\rho_\delta(s) = \frac{1}{\delta^4}\rho\left(\frac{|s|}{\delta}\right), \qquad \int_0^1 \rho(s)s^3 ds = \frac{2}{\pi}.$$

Hereafter we denote by $a = \delta\sqrt{m^2 + n^2}$. Noting that $\lambda_0(m,n) = -m^2 - n^2$, the above inequality yields $0 \leq |\lambda_\delta(m,n)| \leq |\lambda_0(m,n)|$.

Since $1 - \cos\theta \geq \frac{\theta^2}{2} - \frac{\theta^4}{24}$ holds for any $\theta \in \mathbb{R}$, we have

$$\delta^2 |\lambda_\delta(m,n)| = \int_0^1 r\rho(r) \int_0^{2\pi} (1 - \cos(ar\cos\theta))\, d\theta dr$$
$$\geq \int_0^1 r\rho(r) \int_0^{2\pi} \left(\frac{a^2 r^2 \cos^2\theta}{2}\right) d\theta dr - \int_0^1 r\rho(r) \int_0^{2\pi} \left(\frac{a^4 r^4 \cos^4\theta}{24}\right) d\theta dr$$
$$\geq a^2 - \frac{a^4}{12} \int_0^1 r\rho(r) \int_0^{2\pi} \frac{r^2 \cos^2\theta}{2} d\theta dr = a^2 - \frac{a^4}{12},$$

where we used the second moment condition twice respectively in the second inequality and the last equality and the fact $r^2 \cos^2\theta \leq 1$ for $\forall r \in [0, 1]$.

*Case A*: $a = \delta\sqrt{m^2 + n^2} \leq \pi$. We obtain

$$S_{mn} = \frac{1}{\delta^2 |\lambda_\delta(m,n)|} - \frac{1}{\delta^2(m^2 + n^2)} \leq \frac{1}{a^2 - \frac{a^4}{12}} - \frac{1}{a^2} = \frac{1}{12 - a^2} \leq \frac{1}{12 - \pi^2}.$$

*Case B*: $a = \delta\sqrt{m^2 + n^2} > \pi$. Denote by

$$I(r;a) = 4 \int_0^{\frac{\pi}{2}} (1 - \cos(ar\cos\theta))\, d\theta, \quad \forall r \in (0, 1).$$

With the fact that $1 - \cos t \geq t^2/4$ for any $t \in [0, 11/4]$, it holds, for $r \in \left(0, \frac{2}{a}\right]$, that

$$I(r;a) \geq \int_0^{\frac{\pi}{2}} a^2 r^2 \cos^2\theta d\theta = \frac{\pi}{4} a^2 r^2 \geq \frac{\pi^3}{4} r^2.$$

For $r \in \left(\frac{2}{a}, 1\right)$, if we denote by $\xi = ar$, then function

$$J(\xi) = I(r;a) = 4 \int_0^{\frac{\pi}{2}} (1 - \cos(\xi \cos\theta))\, d\theta$$

is increasing for $\xi \in (2, \pi)$, thus

$$I(a;r) = J(\xi) \geq J(2) \geq J(2) r^2.$$

For $\xi \in (\pi, a)$, i.e., $r \in \left(\frac{\pi}{a}, 1\right)$,

$$I(r;a) = J(\xi) = \frac{4}{\xi} \int_0^\xi \frac{1 - \cos t}{1 - (t/\xi)^2} dt \geq \frac{4}{\xi} \int_0^\xi (1 - \cos t) dt$$
$$= \frac{4}{\xi}(\xi - \sin\xi) \geq 4 - \frac{4}{\xi} > 2r^2.$$

By taking $C = 1$ such that $C \leq \min\left\{\frac{\pi^3}{4}, J(2), 2\right\} \approx 6/5$, we then get

$$I(r;a) \geq r^2$$

for any $r \in (0,1)$ and $a \geq \pi$. Thus,

$$\delta^2 |\lambda_\delta(m,n)| = 4 \int_0^1 \int_0^{\frac{\pi}{2}} r\rho(r) \left(1 - \cos(\delta r \sqrt{m^2 + n^2} \cos\theta)\right) d\theta dr$$

$$\geq 4 \int_0^1 r^3 \rho(r) dr = \frac{8}{\pi}.$$

As a result

$$S_{mn} = \frac{1}{\delta^2 |\lambda_\delta(m,n)|} - \frac{1}{\delta^2(m^2+n^2)} \leq \frac{1}{\delta^2 |\lambda_\delta(m,n)|} \leq \frac{\pi}{8}.$$

Hence we obtain the desired result by considering the above two cases. □

### 5.5. *Proof of Theorem 4.1*

**Proof.** We denote

$$E(t) = \|\phi_N^\delta(t,\cdot) - \phi_N^0(t,\cdot)\|.$$

Then it is easy to verify

$$\frac{1}{2}\frac{d}{dt}E^2 = \left(\frac{\partial}{\partial t}\phi_N^\delta - \frac{\partial}{\partial t}\phi_N^0, \phi_N^\delta - \phi_N^0\right).$$

Substituting (3.6) into the above equation yields that

$$E\frac{dE}{dt} = \left(\tilde{\mathcal{L}}_\delta \mathcal{G}_\delta \phi_N^\delta - \tilde{\mathcal{L}}_0 \mathcal{G}_0 \phi_N^0, \phi_N^\delta - \phi_N^0\right) + \left(P_N[\tilde{\mathcal{L}}_\delta F'(\phi_N^\delta)] - P_N[\tilde{\mathcal{L}}_0 F'(\phi_N^0)], \phi_N^\delta - \phi_N^0\right)$$
$$:= I + II,$$

where

$$I = \left(\tilde{\mathcal{L}}_\delta \mathcal{G}_\delta \phi_N^\delta - \tilde{\mathcal{L}}_0 \mathcal{G}_0 \phi_N^0, \phi_N^\delta - \phi_N^0\right)$$
$$= \left((\tilde{\mathcal{L}}_\delta \mathcal{G}_\delta - \tilde{\mathcal{L}}_0 \mathcal{G}_0)\phi_N^0, \phi_N^\delta - \phi_N^0\right) + \left(\tilde{\mathcal{L}}_\delta \mathcal{G}_\delta(\phi_N^\delta - \phi_N^0), \phi_N^\delta - \phi_N^0\right)$$
$$\leq \|(\tilde{\mathcal{L}}_\delta \mathcal{G}_\delta - \tilde{\mathcal{L}}_0 \mathcal{G}_0)\phi_N^0\| E + \left(\tilde{\mathcal{L}}_\delta \mathcal{G}_\delta(\phi_N^\delta - \phi_N^0), \phi_N^\delta - \phi_N^0\right)$$

and

$$II = \left(\tilde{\mathcal{L}}_\delta F'(\phi_N^\delta) - \tilde{\mathcal{L}}_0 F'(\phi_N^0), \phi_N^\delta - \phi_N^0\right)$$
$$= \left(\tilde{\mathcal{L}}_\delta F'(\phi_N^\delta) - \mathcal{L}_\delta F'(\phi_N^0), \phi_N^\delta - \phi_N^0\right) + \left(\tilde{\mathcal{L}}_\delta F'(\phi_N^0) - \tilde{\mathcal{L}}_0 F'(\phi_N^0), \phi_N^\delta - \phi_N^0\right)$$
$$:= II_1 + II_2.$$

Note that

$$II_1 = \left(\tilde{\mathcal{L}}_\delta F'(\phi_N^\delta) - \tilde{\mathcal{L}}_\delta F'(\phi_N^0), \phi_N^\delta - \phi_N^0\right)$$
$$= \left(F'(\phi_N^\delta) - F'(\phi_N^0), \tilde{\mathcal{L}}_\delta(\phi_N^\delta - \phi_N^0)\right)$$
$$\leq \frac{1}{2}\|F'(\phi_N^\delta) - F'(\phi_N^0)\|^2 + \frac{1}{2}\|\tilde{\mathcal{L}}_\delta(\phi_N^\delta - \phi_N^0)\|^2$$
$$\leq \frac{\kappa^2}{2}E^2 + \frac{1}{2}\left(\tilde{\mathcal{L}}_\delta^2(\phi_N^\delta - \phi_N^0), \phi_N^\delta - \phi_N^0\right)$$

where we used Young's inequality in the first inequality and Lipschitz continuity of $F'(\phi)$ in the second inequality, and

$$II_2 = \left(\tilde{\mathcal{L}}_\delta F'(\phi_N^0) - \tilde{\mathcal{L}}_0 F'(\phi_N^0), \phi_N^\delta - \phi_N^0\right)$$
$$\leq \|(\tilde{\mathcal{L}}_\delta - \tilde{\mathcal{L}}_0)F'(\phi_N^0)\|E = \|\mathcal{L}_0^{-2}(\tilde{\mathcal{L}}_\delta - \tilde{\mathcal{L}}_0)\mathcal{L}_0^2 F'(\phi_N^0)\|E$$
$$\leq C\delta^2 \|\mathcal{L}_0^2 F'(\phi_N^0)\|E.$$

Thus we have

$$E\frac{dE}{dt} \leq \frac{\kappa^2}{2}E^2 + \|(\tilde{\mathcal{L}}_\delta\mathcal{G}_\delta - \tilde{\mathcal{L}}_0\mathcal{G}_0)\phi_N^0\|E + C\delta^2\|\mathcal{L}_0 F'(\phi_N^0)\|E$$
$$+ \left(\tilde{\mathcal{L}}_\delta\left(\mathcal{G}_\delta + \frac{1}{2}\tilde{\mathcal{L}}_\delta\right)(\phi_N^\delta - \phi_N^0), \phi_N^\delta - \phi_N^0\right).$$

Specifically with $\tilde{\mathcal{L}}_\delta = \mathcal{L}_\delta$ and $\tilde{\mathcal{L}}_0 = \mathcal{L}_0$, by referring to expressions of $\mathcal{G}_\delta$ and $\mathcal{G}_0$, we have

$$\|(\tilde{\mathcal{L}}_\delta\mathcal{G}_\delta - \tilde{\mathcal{L}}_0\mathcal{G}_0)\phi_N^0\| \leq C\delta^2(\|\mathcal{L}_0^4\phi_N^0\| + \|\mathcal{L}_0^3\phi_N^0\| + \|\mathcal{L}_0^2\phi_N^0\|).$$

Then it yields that

$$E\frac{dE}{dt} \leq \frac{\kappa^2}{2}E^2 + C\delta^2(\|\mathcal{L}_0^4\phi_N^0\| + \|\mathcal{L}_0^3\phi_N^0\| + \|\mathcal{L}_0^2\phi_N^0\| + \|\mathcal{L}_0 F'(\phi_N^0)\|)E$$
$$+ \left(\mathcal{L}_\delta\left(I + \mathcal{L}_\delta^2 + \frac{5}{2}\mathcal{L}_\delta\right)(\phi_N^\delta - \phi_N^0), \phi_N^\delta - \phi_N^0\right).$$

To estimate $g(\lambda_\delta) = \frac{5}{2}\lambda_\delta^2 + \lambda_\delta^3 + \lambda_\delta$ with $\lambda_\delta < 0$, we note that $g(\lambda_\delta) \leq 0$ when $\lambda_\delta \leq -\frac{5}{2}$. This yields that

$$\frac{dE}{dt} \leq \frac{\kappa^2}{2}E + C\delta^2(\|\mathcal{L}_0^4\phi_N^0\| + \|\mathcal{L}_0^3\phi_N^0\| + \|\mathcal{L}_0^2\phi_N^0\| + \|\mathcal{L}_0 F'(\phi_N^0)\|).$$

For the case $-\frac{5}{2} < \lambda_\delta < 0$, we note that $g(-\frac{5}{2}) < 0$ is finite and independent of $\delta$, $g(0) = 0$. Due to the fact that $g(\lambda_\delta)$ is continuous on the interval $[-\frac{5}{2}, 0]$, we know that $g(\lambda_\delta)$ is bounded for $-\frac{5}{2} < \lambda_\delta < 0$. Thus we have

$$\frac{dE}{dt} \leq \left(\frac{\kappa^2}{2} + C\right)E + C\delta^2(\|\mathcal{L}_0^4\phi_N^0\| + \|\mathcal{L}_0^3\phi_N^0\| + \|\mathcal{L}_0^2\phi_N^0\| + \|\mathcal{L}_0 F'(\phi_N^0)\|).$$

Again Gronwall's inequality yields the desired estimate. □

### 5.6. *Proof of Theorem 4.2*

**Proof.** With errors $e_\phi^n = \phi_N^{\delta,n} - \phi_N^{0,n}$, $e_\mu^n = \mu_N^{\delta,n} - \mu_N^{0,n}$, $e_r^n = r^{\delta,n} - r^{0,n}$, error equations are written as

$$\left(3e_\phi^{n+1} - 4e_\phi^n + e_\phi^{n-1}, q\right) - 2\tau\left(\tilde{\mathcal{L}}_\delta e_\mu^{n+1}, q\right) = 2\tau\left((\tilde{\mathcal{L}}_\delta - \mathcal{L}_0)\mu_N^{0,n+1}, q\right),$$
$$\left(e_\mu^{n+1}, w\right) = \left(\mathcal{G}_\delta e_\phi^{n+1}, w\right) + \left((\mathcal{G}_\delta - \mathcal{G}_0)\phi_N^{0,n+1}, w\right) + \beta(e_\phi^{n+1}, w)$$

$$+ \left(e_r^{n+1} H(\bar{\phi}_N^{\delta,n+1}), w\right) + r^{0,n+1}\left(H(\bar{\phi}_N^{\delta,n+1}) - H(\bar{\phi}_N^{0,n+1}), w\right),$$

$$3e_r^{n+1} - 4e_r^n + e_r^{n-1} = \frac{1}{2}\left(H(\bar{\phi}_N^{\delta,n+1}), 3e_\phi^{n+1} - 4e_\phi^n + e_\phi^{n-1}\right)$$
$$+ \frac{1}{2}\left(H(\bar{\phi}_N^{\delta,n+1}) - H(\bar{\phi}_N^{0,n+1}), 3\phi_N^{0,n+1} - 4\phi_N^{0,n} + \phi_N^{0,n-1}\right).$$

Testing with $q = e_\mu^{n+1}$ in the first equation, $w = 3e_\phi^{n+1} - 4e_\phi^n + e_\phi^{n-1}$ in the second equation and multiplying the third equation with $2e_r^{n+1}$ lead to

$$\left(3e_\phi^{n+1} - 4e_\phi^n + e_\phi^{n-1}, e_\mu^{n+1}\right) + 2\tau\|(-\tilde{\mathcal{L}}_\delta)^{\frac{1}{2}} e_\mu^{n+1}\|^2 = 2\tau\left((\tilde{\mathcal{L}}_\delta - \mathcal{L}_0)\mu_N^{0,n+1}, e_\mu^{n+1}\right),$$

$$\left(e_\mu^{n+1}, 3e_\phi^{n+1} - 4e_\phi^n + e_\phi^{n-1}\right) = \left(\mathcal{G}_\delta e_\phi^{n+1}, 3e_\phi^{n+1} - 4e_\phi^n + e_\phi^{n-1}\right)$$
$$+ \left((\mathcal{G}_\delta - \mathcal{G}_0)\phi_N^{0,n+1}, 3e_\phi^{n+1} - 4e_\phi^n + e_\phi^{n-1}\right)$$
$$+ \beta(e_\phi^{n+1}, 3e_\phi^{n+1} - 4e_\phi^n + e_\phi^{n-1})$$
$$+ \left(e_r^{n+1} H(\bar{\phi}_N^{\delta,n+1}), 3e_\phi^{n+1} - 4e_\phi^n + e_\phi^{n-1}\right)$$
$$+ r^{0,n+1}\left(H(\bar{\phi}_N^{\delta,n+1}) - H(\bar{\phi}_N^{0,n+1}), 3e_\phi^{n+1} - 4e_\phi^n + e_\phi^{n-1}\right),$$

$$2e_r^{n+1}(3e_r^{n+1} - 4e_r^n + e_r^{n-1}) = e_r^{n+1}\left(H(\bar{\phi}_N^{\delta,n+1}), 3e_\phi^{n+1} - 4e_\phi^n + e_\phi^{n-1}\right)$$
$$+ e_r^{n+1}\left(H(\bar{\phi}_N^{\delta,n+1}) - H(\bar{\phi}_N^{0,n+1}), 3\phi_N^{0,n+1} - 4\phi_N^{0,n} + \phi_N^{0,n-1}\right).$$

Summing up the above equations yields

$$2\tau\|(-\tilde{\mathcal{L}}_\delta)^{\frac{1}{2}} e_\mu^{n+1}\|^2 + \frac{1}{2}\Big(\|(I+\mathcal{L}_\delta)e_\phi^{n+1}\|^2 - \|(I+\mathcal{L}_\delta)e_\phi^n\|^2 + \|(I+\mathcal{L}_\delta)(2e_\phi^{n+1} - e_\phi^n)\|^2$$
$$+ \|(I+\mathcal{L}_\delta)(e_\phi^{n+1} - 2e_\phi^n + e_\phi^{n-1})\|^2 - \|(I+\mathcal{L}_\delta)(2e_\phi^n - e_\phi^{n-1})\|^2\Big)$$
$$+ \frac{\beta}{2}\left(\|e_\phi^{n+1}\|^2 - \|e_\phi^n\|^2 + \|2e_\phi^{n+1} - e_\phi^n\|^2 + \|e_\phi^{n+1} - 2e_\phi^n + e_\phi^{n-1}\|^2 - \|2e_\phi^n - e_\phi^{n-1}\|^2\right)$$
$$+ |e_r^{n+1}|^2 - |e_r^n|^2 + |2e_r^{n+1} - e_r^n|^2 + |e_r^{n+1} - 2e_r^n + e_r^{n-1}|^2 - |2e_r^n - e_r^{n-1}|^2$$
$$= \left\{2\tau\left((\tilde{\mathcal{L}}_\delta - \mathcal{L}_0)\mu_N^{0,n+1}, e_\mu^{n+1}\right)\right\} + \left\{-\left((\mathcal{G}_\delta - \mathcal{G}_0)\phi_N^{0,n+1}, 3e_\phi^{n+1} - 4e_\phi^n + e_\phi^{n-1}\right)\right\}$$
$$+ \left\{-r^{0,n+1}\left(H(\bar{\phi}_N^{\delta,n+1}) - H(\bar{\phi}_N^{0,n+1}), 3e_\phi^{n+1} - 4e_\phi^n + e_\phi^{n-1}\right)\right\}$$
$$+ \left\{e_r^{n+1}\left(H(\bar{\phi}_N^{\delta,n+1}) - H(\bar{\phi}_N^{0,n+1}), 3\phi_N^{0,n+1} - 4\phi_N^{0,n} + \phi_N^{0,n-1}\right)\right\}$$
$$:= \sum_{i=1}^4 I_i + 2\tau(\mathcal{L}_\delta e_\mu^{n+1}, e_\mu^{n+1}).$$

We next estimate the four terms one by one. By $\tilde{\mathcal{L}}_\delta = \mathcal{L}_\delta$, the first term can be controlled as

$$|I_1| = |2\tau\left((\mathcal{L}_\delta - \mathcal{L}_0)\mu_N^{0,n+1}, e_\mu^{n+1}\right)|$$
$$= \left|2\tau\left((\mathcal{L}_\delta - \mathcal{L}_0)(-\mathcal{L}_\delta)^{-\frac{1}{2}}\mu_N^{0,n+1}, (-\mathcal{L}_\delta)^{\frac{1}{2}} e_\mu^{n+1}\right)\right|$$

$$\leq C\tau\delta^4 + \frac{\tau}{6}\|(-\mathcal{L}_\delta)^{\frac{1}{2}}e_\mu^{n+1}\|^2$$

where the inequality holds due to the Lemma 4.1. The second term is bounded by

$$|I_2| = \left|\left((\mathcal{G}_\delta - \mathcal{G}_0)\phi_N^{0,n+1}, 3e_\phi^{n+1} - 4e_\phi^n + e_\phi^{n-1}\right)\right|$$
$$\leq \left|2\tau\left((\mathcal{G}_\delta - \mathcal{G}_0)\phi_N^{0,n+1}, \mathcal{L}_\delta e_\mu^{n+1}\right)\right| + \left|2\tau\left((\mathcal{G}_\delta - \mathcal{G}_0)\phi_N^{0,n+1}, (\tilde{\mathcal{L}}_\delta - \mathcal{L}_0)\mu_N^{0,n+1}\right)\right|$$
$$= \left|2\tau\left((\mathcal{G}_\delta - \mathcal{G}_0)(-\tilde{\mathcal{L}}_\delta)^{\frac{1}{2}}\phi_N^{0,n+1}, (-\tilde{\mathcal{L}}_\delta)^{\frac{1}{2}}e_\mu^{n+1}\right)\right|$$
$$+ \left|2\tau\left((\mathcal{G}_\delta - \mathcal{G}_0)\phi_N^{0,n+1}, (\tilde{\mathcal{L}}_\delta - \mathcal{L}_0)\mu_N^{0,n+1}\right)\right|$$
$$\leq C\tau\delta^4 + \frac{\tau}{6}\|(-\tilde{\mathcal{L}}_\delta)^{\frac{1}{2}}e_\mu^{n+1}\|^2$$

where we use the expression $\mathcal{G}_\delta - \mathcal{G}_0 = (\mathcal{L}_\delta - \mathcal{L}_0)(\mathcal{L}_\delta + \mathcal{L}_0 + 2I)$. Similarly, the third term can be estimated as

$$|I_3| = \left|r^{0,n+1}\left(H(\bar{\phi}_N^{\delta,n+1}) - H(\bar{\phi}_N^{0,n+1}), 3e_\phi^{n+1} - 4e_\phi^n + e_\phi^{n-1}\right)\right|$$
$$\leq \left|2r^{0,n+1}\tau\left(H(\bar{\phi}_N^{\delta,n+1}) - H(\bar{\phi}_N^{0,n+1}), \tilde{\mathcal{L}}_\delta e_\mu^{n+1}\right)\right|$$
$$+ \left|2r^{0,n+1}\tau\left(H(\bar{\phi}_N^{\delta,n+1}) - H(\bar{\phi}_N^{0,n+1}), (\tilde{\mathcal{L}}_\delta - \mathcal{L}_0)\mu_N^{0,n+1}\right)\right|$$
$$= \left|2r^{0,n+1}\tau\left((-\tilde{\mathcal{L}}_\delta)^{\frac{1}{2}}(H(\bar{\phi}_N^{\delta,n+1}) - H(\bar{\phi}_N^{0,n+1})), (-\tilde{\mathcal{L}}_\delta)^{\frac{1}{2}}e_\mu^{n+1}\right)\right|$$
$$+ \left|2r^{0,n+1}\tau\left(H(\bar{\phi}_N^{\delta,n+1}) - H(\bar{\phi}_N^{0,n+1}), (\tilde{\mathcal{L}}_\delta - \mathcal{L}_0)\mu_N^{0,n+1}\right)\right|$$
$$\leq C\tau(\|H(\bar{\phi}_N^{\delta,n+1}) - H(\bar{\phi}_N^{0,n+1})\|^2 + \|(-\tilde{\mathcal{L}}_\delta)^{\frac{1}{2}}(H(\bar{\phi}_N^{\delta,n+1}) - H(\bar{\phi}_N^{0,n+1}))\|^2)$$
$$+ \frac{\tau}{6}\|(-\tilde{\mathcal{L}}_\delta)^{\frac{1}{2}}e_\mu^{n+1}\|^2 + C\tau\delta^4$$
$$\leq C\tau(\|2e_\phi^n - e_\phi^{n-1}\|^2 + \|(I + \mathcal{L}_\delta)(2e_\phi^n - e_\phi^{n-1})\|^2) + \frac{\tau}{6}\|\mathcal{L}_\delta e_\mu^{n+1}\|^2 + C\tau\delta^4,$$

where in the last second inequality we note with

$$\mathcal{W}(\phi) = \int_\Omega \{F(\phi)\} - \frac{\beta}{2}\phi^2 + C_H$$

that

$$H(\bar{\phi}_N^{\delta,n+1}) - H(\bar{\phi}_N^{0,n+1}) = \frac{U(\bar{\phi}_N^{\delta,n+1})}{\sqrt{\mathcal{W}(\bar{\phi}_N^{\delta,n+1})}} - \frac{U(\bar{\phi}_N^{0,n+1})}{\sqrt{\mathcal{W}(\bar{\phi}_N^{0,n+1})}}$$
$$= \frac{U(\bar{\phi}_N^{\delta,n+1}) - U(\bar{\phi}_N^{0,n+1})}{\sqrt{\mathcal{W}(\bar{\phi}_N^{\delta,n+1})}} + \frac{U(\bar{\phi}_N^{0,n+1})}{\sqrt{\mathcal{W}(\bar{\phi}_N^{\delta,n+1})}\sqrt{\mathcal{W}(\bar{\phi}_N^{0,n+1})}} \frac{\mathcal{W}(\bar{\phi}_N^{\delta,n+1}) - \mathcal{W}(\bar{\phi}_N^{0,n+1})}{\sqrt{\mathcal{W}(\bar{\phi}_N^{\delta,n+1})} + \sqrt{\mathcal{W}(\bar{\phi}_N^{0,n+1})}}.$$

Since $\left|\sqrt{\mathcal{W}(\bar{\phi}_N^{\delta,n+1})}\right|, \left|\sqrt{\mathcal{W}(\bar{\phi}_N^{0,n+1})}\right|, |U(\bar{\phi}_N^{0,n+1})|$ are bounded, using Lemma 3.2 with above splitting yields

$$\|(-\mathcal{L}_\delta)^{\frac{1}{2}}(H(\bar{\phi}_N^{\delta,n+1}) - H(\bar{\phi}_N^{0,n+1}))\|^2$$

$$=(-\mathcal{L}_\delta(H(\bar{\phi}_N^{\delta,n+1}) - H(\bar{\phi}_N^{0,n+1})), H(\bar{\phi}_N^{\delta,n+1}) - H(\bar{\phi}_N^{0,n+1}))$$
$$\leq C(-\mathcal{L}_\delta(U(\bar{\phi}_N^{\delta,n+1}) - U(\bar{\phi}_N^{0,n+1})), U(\bar{\phi}_N^{\delta,n+1}) - U(\bar{\phi}_N^{0,n+1}))$$
$$+C(-\mathcal{L}_\delta(\mathcal{W}(\bar{\phi}_N^{\delta,n+1}) - \mathcal{W}(\bar{\phi}_N^{0,n+1})), \mathcal{W}(\bar{\phi}_N^{\delta,n+1}) - \mathcal{W}(\bar{\phi}_N^{0,n+1}))$$
$$\leq C(-\tilde{\mathcal{L}}_\delta(2e_\phi^n - e_\phi^{n-1}), 2e_\phi^n - e_\phi^{n-1}) + C(2e_\phi^n - e_\phi^{n-1}, 2e_\phi^n - e_\phi^{n-1})$$
$$\leq C(\|(I+\mathcal{L}_\delta)(2e_\phi^n - e_\phi^{n-1}))\|^2 + \|2e_\phi^n - e_\phi^{n-1}\|^2)$$

where the last inequality follows from the inequality (5.4). The last term $I_4$ is controlled by

$$|I_4| = \left| e_r^{n+1} \left( H(\bar{\phi}_N^{\delta,n+1}) - H(\bar{\phi}_N^{0,n+1}), 3\phi_N^{0,n+1} - 4\phi_N^{0,n} + \phi_N^{0,n-1} \right) \right|$$
$$= \left| 2\tau \left( H(\bar{\phi}_N^{\delta,n+1}) - H(\bar{\phi}_N^{0,n+1}), e_r^{n+1} \mathcal{L}_0 \mu_N^{0,n+1} \right) \right|$$
$$\leq C\tau \|2e_\phi^n - e_\phi^{n-1}\|^2 + C\tau |e_r^{n+1}|^2.$$

Then we arrive at

$$2\tau\|(-\mathcal{L}_\delta)^{\frac{1}{2}} e_\mu^{n+1}\|^2 + \frac{1}{2}\|(I+\mathcal{L}_\delta)e_\phi^{n+1}\|^2 + \left(\frac{1}{2} - C\tau\right) \|(I+\mathcal{L}_\delta)(2e_\phi^{n+1} - e_\phi^n)\|^2$$
$$+ (1 - C\tau)|e_r^{n+1}|^2 + \left(\frac{\beta}{2} - C\tau\right)(\|e_\phi^{n+1}\|^2 + \|2e_\phi^{n+1} - e_\phi^n\|^2)$$
$$\leq \frac{\tau}{2}\|(-\mathcal{L}_\delta)^{\frac{1}{2}} e_\mu^{n+1}\|^2 + \frac{1}{2}\|(I+\mathcal{L}_\delta)e_\phi^n\|^2 + \left(\frac{1}{2} + C\tau\right) \|(I+\mathcal{L}_\delta)(2e_\phi^n - e_\phi^{n-1})\|^2$$
$$+ (1 + C\tau)|e_r^n|^2 + \left(\frac{\beta}{2} + C\tau\right)(\|e_\phi^n\|^2 + \|2e_\phi^n - e_\phi^{n-1}\|^2) + C\tau\delta^4.$$

Finally applying the discrete Gronwall's inequality on the above inequality leads to the desired estimate

$$\|e_\phi^{n+1}\| + |e_r^{n+1}| + \|(I+\mathcal{L}_\delta)e_\phi^{n+1}\| \leq C(T,\phi_0)\delta^2. \qquad \square$$